\def\Datum{July 30, 2005}
\magnification=\magstephalf
\hsize=16true cm
\advance\vsize5truemm
\parskip4pt plus 1pt
\frenchspacing

\vglue6truemm

\parindent25pt
\newif\ifdraft\newif\ifneu
\newdimen\SkIp\SkIp=\parskip
\pageno=1

\font\tenrsfs=rsfs10
\font\sevenrsfs=rsfs10 scaled 700
\font\fiversfs=rsfs10 scaled 500
\font\tenbf=txb

\newfam\rsfsfam
\textfont\rsfsfam=\tenrsfs
\scriptfont\rsfsfam=\sevenrsfs
\scriptscriptfont\rsfsfam=\fiversfs
\def\rsfs{\fam\rsfsfam}

\def\Ues#1{{\medskip\noindent\bf#1.}}
\def\beginabs{\begingroup\parindent0pt\openup-1pt}
\long\def\endabs{\medskip\endgroup}

\def\GaRGaR{1}
\ifx\ttex\GaRGaR 
 \font\eightrm=txr scaled 800
 \font\eightsl=txsl scaled 800
 \font\eightbf=txb scaled 800
 \font\ninerm=txr scaled 900
 \font\ninesl=txsl scaled 900
 \font\ninebf=txb scaled 900
 
 \font\gross=txb scaled\magstep1
 \font\Gross=txb scaled\magstep2
\else 
 \font\eightrm=cmr8
 \font\eightsl=cmsl8
 \font\eightbf=cmbx8
 \font\ninerm=cmr9
 \font\ninesl=cmsl9
 \font\ninebf=cmbx9
 
 \font\gross=cmbx10 scaled\magstep1
 \font\Gross=cmbx10 scaled\magstep2
\fi

\def\nline{\hfill\break}
\long\def\fussnote#1#2{{\baselineskip=9pt
     \setbox\strutbox=\hbox{\vrule height 7 pt depth 2pt width 0pt}%
     \eightrm
     \footnote{#1}{#2}}}

\def\footnoterule{\kern-3pt
         \hrule width 2 true cm
         \kern 2.6pt}

\def\Times{\]2\times\]2}

\def\mapright#1{\smash{\mathop{\hbox to 35pt{\rightarrowfill}}\limits^{#1}}}

\def\cd{{\cdot}}
\font\trm=cmr17
\font\eighti=cmmi8
\font\nineti=cmmi9
\font\Sans=cmss10
\font\Icke=cmcsc10
\font\eightt=cmtt8
\font\eightsy=cmsy8
\font\ninesy=cmsy9

\font\tenmsb=msbm10\font\sevenmsb=msbm7\font\fivemsb=msbm5
\newfam\msbfam
\textfont\msbfam=\tenmsb\scriptfont\msbfam=\sevenmsb
\scriptscriptfont\msbfam=\fivemsb
\def\msb{\fam\msbfam}

\def\ltimes{{\msb\char110}}

\font\tenmib=cmmib10 \font\eightmib=cmmib8 \font\fivemib=cmmib5
\newfam\mibfam
\textfont\mibfam=\tenmib \scriptfont\mibfam=\eightmib
\scriptscriptfont\mibfam=\fivemib
\def\mib{\fam\mibfam\tenmib}

\font\tenmi=cmmi10 \font\eightmi=cmmi8 \font\fivemi=cmmi5
\newfam\mifam
\textfont\mifam=\tenmi \scriptfont\mifam=\eightmi
\scriptscriptfont\mifam=\fivemi

\font\tenfr=eufm10 \font\eightfr=eufm8 \font\fivefr=eufm5
\newfam\frfam
\textfont\frfam=\tenfr \scriptfont\frfam=\eightfr
\scriptscriptfont\frfam=\fivefr
\def\frak{\fam\frfam\tenfr}

\font\tenbfr=eufb10 \font\eightbfr=eufb8 \font\fivebfr=eufb5
\newfam\bfrfam
\textfont\bfrfam=\tenbfr \scriptfont\bfrfam=\eightbfr
\scriptscriptfont\bfrfam=\fivebfr
\def\bfrak{\fam\bfrfam\tenbfr}

\def\Quot#1#2{\raise 2pt\hbox{$#1\mskip-1.2\thinmuskip$}\big/%
     \lower2pt\hbox{$\mskip-0.8\thinmuskip#2$}}
\def\ddz#1{\raise1.5pt\hbox{$\[1\partial\!$}/\raise-2.5pt\hbox{$\!\partial z_{#1}$}}
\def\ddw{\raise1.5pt\hbox{$\[1\partial\!$}/\raise-2.5pt\hbox{$\!\partial w\,$}}
\def\ddo{\raise1.5pt\hbox{$\,\partial\!$}/\raise-2.5pt\hbox{$\!\partial \overline z\,$}}
\def\hol{\7{hol}}\def\aut{\7{aut}}

\def\Klein{\ninerm\textfont1=\nineti\textfont0=\ninerm\def\sl{\ninesl}
\textfont2=\ninesy\baselineskip11.2pt\def\bf{\ninebf}}
\def\klein{\eightrm\textfont1=\eighti\textfont0=\eightrm\def\sl{\eightsl}
\textfont2=\eightsy\baselineskip10.5pt\def\bf{\eightbf}}
\def\Matrix#1#2#3#4{{\klein\Big(\matrix{#1&\!\!#2\cr\noalign{\vskip-2pt}#3&\!\!#4\cr}\Big)}}

\def\Kl#1{\raise1pt\hbox{$\scriptstyle($}#1\raise1pt\hbox{$\scriptstyle)$}}

\newcount\ite\ite=1
\def\0{\global\ite=1\1}
\def\1{\item{\rm(\romannumeral\the\ite)}\2}
\def\2{\global\advance\ite1}
\def\3#1{{\mib#1}}
\def\5#1{{\cal#1}}
\def\6#1{{\rsfs#1}}
\def\7#1{\mathop{\frak#1\[2}\nolimits}
\def\8#1#2#3#4{\big(\raise-2pt\hbox{\rlap{\raise6pt%
  \hbox{$\scriptstyle#1\[1#2$}}\hbox{$\scriptstyle#3\[2#4$}}\big)}\
\def\9#1{\rlap{\hskip.2pt$#1$}\rlap{\hskip.4pt$#1$}#1}

\def\[#1{\mskip#1mu}
\def\]#1{\mskip-#1mu}

\def\tilde{\widetilde}

\def\<{\;\;\Longleftrightarrow\;\;}
\def\To#1#2{(\romannumeral#1) $\Longrightarrow$ (\romannumeral#2)}

\def\epsilon{\varepsilon}
\def\hat{\rlap{\raise-1.5pt\hbox{\hskip2.5pt\trm\char94}}}
\def\phi{\varphi}
\def\steil#1{\hbox{\rm~~#1~~}}
\def\Steil#1{\hbox{\rm\quad #1\quad}}
\def\p{p.\nobreak\hskip2pt}
\def\Im{\hbox{\rm Im}\[1}
\def\Re{\hbox{\rm Re}\[1}

\def\MOP#1{\expandafter\edef\csname #1\endcsname{%
   \mathop{\hbox{\Sans #1}}\nolimits
    }}

\MOP{Aut}\MOP{ad}\MOP{GL}\MOP{SU}\MOP{SO}\MOP{tr}\MOP{Sp}\MOP{Spin}
\MOP{SL}\MOP{id}\MOP{Ad}\MOP{R}\MOP{PSL}\MOP{U}\MOP{O}\MOP{Fix}

\def\sqr#1#2{{\,\vcenter{\vbox{\hrule height.#2pt\hbox{\vrule width.#2pt
height#1pt \kern#1pt\vrule width.#2pt}\hrule height.#2pt}}\,}}
\def\qed{\hfill\ifmmode\sqr66\else$\sqr66$\par\fi\rm}

\def\One{{1\kern-3.8pt 1}}
\def\one{{1\kern-3.1pt 1}}

\font\tenmia=txmia
\font\sevenmia=txmia scaled 700
\font\fivemia=txmia scaled 500

\def\Bbord#1{\mathchoice
{\hbox{\tenmia\char#1}}
{\hbox{\tenmia\char#1}}
{\hbox{\sevenmia\char#1}}
{\hbox{\fivemia\char#1}}
}

\def\CC{\Bbord{131}}
\def\NN{\Bbord{142}}
\def\PP{\Bbord{144}}
\def\RR{\Bbord{146}}
\def\ZZ{\Bbord{154}}

\newcount\nummer
\newcount\parno\parno=0

\def\KAP#1#2{\par\bigbreak\global\advance\parno by1%
\def\test{#1}\ifx\test\empty\else%
\expandafter\let\csname#1\endcsname=\relax%
\immediate\write\aux{\def\csname#1\endcsname{\the\parno}}%
\expandafter\xdef\csname#1\endcsname{\the\parno}\fi%
{\vskip3pt\noindent\gross\the\parno. #2\hfil}\vskip3pt\rm\nummer=0\nobreak}

\newwrite\aux

\def\Randmark#1{\vadjust{\vbox to 0pt{\vss\hbox to\hsize%
{\fiverm\hskip\hsize\hskip1em\raise 2.5pt\hbox{#1}\hss}}}}

\def\PrN{\the\parno.\the\nummer}

\def\Write#1#2{\global\advance\nummer1\def\test{#1}\ifx\test\empty\else%
\ifdraft\Randmark{#1}\fi\expandafter\let\csname#1\endcsname=\relax%
\immediate\write\aux{\def\csname#1\endcsname{\PrN}}%
\expandafter\xdef\csname#1\endcsname{\PrN}\fi#2}

\def\Num#1{\Write{#1}{\PrN}}
\def\Leqno#1{\Write{#1}{\leqno(\PrN)}}

\def\Proposition#1{{\smallbreak\noindent\bf\Num{#1} Proposition.~}\parskip0pt\sl}
\def\Lemma#1{{\par\noindent\bf\Num{#1} Lemma.~}\parskip0pt\sl}
\def\Corollary#1{{\smallbreak\noindent\bf\Num{#1} Corollary.~}\parskip0pt\sl}

\def\Example#1{{\smallbreak\noindent\bf\Num{#1} Example.~}}

\long\def\Proof{\smallskip\noindent\sl Proof.\parskip\SkIp\rm~~}
\def\Formend{\par\parskip\SkIp\rm}

\def\ruf#1{{{\expandafter\ifx\csname#1\endcsname\relax\xdef\flAG{}%
\message{*** #1 nicht definiert!! ***}\ifdraft\Randmark{#1??}\fi\else%
\xdef\flAG{1}\fi}\ifx\flAG\empty{\bf??}\else\rm\csname#1\endcsname\fi}}
\def\Ruf#1{{\rm(\ruf{#1})}}

\newcount\lit\lit=1
\def\Ref#1{\item{\the\lit.}\expandafter\ifx\csname#1ZZZ\endcsname\relax%
\message{ >>> \the\lit. = #1  <<< }%
\fi%
\expandafter\let\csname#1\endcsname=\relax%
\immediate\write\aux{\def\csname#1\endcsname{\the\lit}}\advance\lit1\ifdraft\Randmark{#1}\fi}

\def\Lit#1{\expandafter\gdef\csname#1ZZZ\endcsname{1}[\ruf{#1}]}
\def\LIT#1#2{\expandafter\gdef\csname#1ZZZ\endcsname{1}[\ruf{#1}, #2]}

\immediate\openout\aux=\jobname.aux

\input epsf


\ifdraft\footline={\hss\sevenrm incomplete draft of \Datum\hss}
\else\nopagenumbers\fi

\headline={\ifnum\pageno>1\sevenrm\ifodd\pageno CR-manifolds of dimension 5
\hss{\tenbf\folio}\else{\tenbf\folio}\hss{Fels-Kaup}\fi\else\hss\fi}

\centerline{\Gross CR-manifolds of dimension 5: A Lie algebra approach}
\bigskip\bigskip
\centerline{Gregor Fels\quad and\quad Wilhelm Kaup} 
{\parindent0pt\footnote{}{\ninerm 2000 Mathematics Subject Classification:
17C50, 32M15, 32M25, 32V10, 32V25.}}

\bigskip\bigskip

\noindent{\bf Abstract.} We study real-analytic Levi degenerate
hypersurfaces $M$ in complex manifolds of dimension $3$, for which the
CR-automorphism group $\Aut(M)$ is a real Lie group acting
transitively on $M$. We provide large classes of examples for such
$M$, compute the corresponding groups $\Aut(M)$ and determine the
maximal subsets of $M$ that cannot be separated by global continuous
CR-functions. It turns out that all our examples, although partly
arising in different contexts, are locally CR-equivalent to the tube
$\6T=\6C\times i\RR^{3}\subset\CC^{3}$ over the future light cone
$\6C:=\{x\in\RR^{3}:x^2_{1}+x^2_{2}=x^2_{3},\,x_{3}>0\}$ in
$3$-dimensional space-time.

\bigskip

\KAP{Introduction}{Introduction} 

The notion of a Cauchy-Riemann-manifold (CR-manifold for short)
generalizes that of a complex manifold. While two complex manifolds of
the same dimension are always locally equivalent, this is no longer
true for CR-manifolds. A basic invariant for a CR-manifold is the
so-called Levi form. In case this form is nondegenerate, {\Icke Chern}
and {\Icke Moser} give in their celebrated paper \Lit{CHMO} further
invariants which characterize up to equivalence the local CR-structure
of real hypersurfaces in $\CC^{n+1}$. For hypersurfaces with {\sl
degenerate} Levi form the corresponding program seems to be much
harder to overcome. Of course, in case the connected real-analytic
(locally closed) hypersurface $M\subset\CC^{n+1}$ has a Levi
nondegenerate point, all points in a dense open subset of $M$ are of
this type. Therefore, of a special interest are those hypersurfaces
which have everywhere degenerate Levi form but are not Levi
flat. Clearly, the smallest $n$ for which hypersurfaces
$M\subset\CC^{n+1}$ with this property can occur is $n=2$, that is,
where $M$ has real dimension $5$. In fact, the tube $\6T: =\6C\times
i\RR^{3}\subset\CC^{3}$ over the future light cone
$\6C:=\{x\in\RR^{3}:x^2_{1}+x^2_{2}=x^2_{3},\,x_{3}>0\}$ in
$3$-dimensional space-time is a well known example of this type. This
CR-manifold, which is the starting point for our studies, is even
homogeneous in the sense that a Lie group of CR-automorphisms acts
transitively. In a way, $\6T$ may be considered as a {\sl model
surface} \Lit{EBEN} for locally homogeneous Levi degenerate surfaces,
similar to the sphere $S^{2n+1}\subset\CC^{n+1}$ for spherical
surfaces in the Levi nondegenerate case, compare e.g. \Lit{BUSH} for
the classification of homogeneous surfaces of this type and \Lit{AZAD}
for general homogeneous CR-manifolds.

In this paper we examine several naturally occurring CR-manifolds of
dimension $5$ with {\sl degenerate} Levi-form. Most (but not all) of
these are locally homogeneous and 2-nondegenerate. Surprisingly, it
turns out that all of the latter are locally equivalent to the tube
$\6T$ over the light cone. From the global point of view, however,
there are large classes of pairwise CR-nonequivalent manifolds which
all are {\sl locally} CR-equivalent to the tube $\6T$. The
classification of all homogeneous CR-manifolds of this type can
essentially be reduced to the study of homogeneous domains in the real
projective space $\PP_3(\RR)$.

{\parindent0pt\medskip The paper is organized as follows:

In Section \ruf{Preliminaries} we fix notation and recall some
basic facts concerning the geometry of CR-manifolds.

In Section \ruf{tube} we recall some (more or less known) properties
of the model example $\6T$ and focus on an explicit description of the
Lie algebra of all infinitesimal CR-automorphisms.

In Section \ruf{local} we investigate hypersurface germs at
$0\in\CC^{3}$ which correspond to the partial normal form (A,ii,2) in
\Lit{EBFT} (the only one out of 8 for holomorphically nondegenerate
hypersurfaces to which uniformly Levi degenerate hypersurfaces can
belong). Besides the fact that a partial normal form presentation of a
hypersurface germ $(M,0)$ is not an invariant in the strict sense
(since it still depends on the choice of an adapted {\sl coordinate
system} and a further group of local CR-automorphisms), in general, it
is rather difficult to deduce from the information encoded in the
partial normal form only whether $(M,0)$ is locally homogeneous or not.
In Proposition \ruf{AR} we give for certain germs a criterion for
local homogeneity. We then introduce a family $(\6M_t)_{t\in \RR}$ of
germs belonging to the partial normal form (A,ii,2) and compute
explicitly the Lie algebras $\hol(\6M_{t},0)$ of of germs of
infinitesimal CR-transformations at $0$. It turns out that
$(\6M_t)_{t\in \RR}$ varies from $\7{sl}(2,\RR)\times\7{su}(2)$
($t<0$) over a non semi-simple Lie algebra ($t=0$) to
$\7{sl}(2,\RR)\times\7{sl}(2,\RR)$ ($t>0$) and jumps for $t=1$ to the
10-dimensional Lie algebra $\7{so}(2,3)$. As main result of the
section we show that actually $\6M_{1}$ is locally CR-equivalent to
the tube $\6T$ and provide an explicit local CR-isomorphism. Moreover,
the family $(\6M_t)_{t\in \RR}$ helps to clarify several results which
can be found in the literature, compare e.g. the concluding paragraph
in Section \ruf{local}.

In Section \ruf{Global} we study global properties of CR-manifolds
that are locally CR-equivalent to the tube $\6T$. Our first result
states that on the universal covering of $\6T$, which is infinitely
sheeted, every global continuous CR-function is the pull-back of a
function defined on $\6T$. We then introduce a CR-deformation family
$(\6T_{t})_{t>0}$ of $\6T$ with interesting properties: The $\6T_{t}$
are pairwise CR-nonequivalent CR-manifolds locally CR-equivalent to
$\6T$ and $\6T_{1}=\6T$. On the other hand, all $\6T_{t}$ are
diffeomorphic to $\6T$ as real manifolds. What concerns the separation
properties by global CR-functions, the $\6T_{t}$ behave differently
for various $t$: For instance, if $t=p/q$ with $p,q\in\NN$ relatively
prime, there are precisely $p$ points in $\6T_{t}$ that cannot be
separated from a given point in $\6T_{t}$ by global CR-functions.
Further, for $t$ irrational, there is a closed hypersurface
$\6N_{t}\subset\6T_{t}$ such that every continuous CR-function on
$\6T_{t}$ is real-analytic outside $\6N_{t}$, while there do exist
continuous CR-functions on $\6T_{t}$ which are not globally
real-analytic.

In Section \ruf{bounded} we consider a hypersurface
$\6R\subset\CC^{3}$ which, in a certain sense, is universal for all
homogeneous CR-manifolds that are locally CR-equivalent to the tube
$\6T$. This hypersurface occurs as the smooth boundary part of the Lie
ball in $\CC^{3}$ (biholomorphic image of Siegel's upper halfplane in
the space of symmetric $2\Times2$-matrices via a Cayley
transformation). One result is that every simply-connected homogeneous
CR-manifold locally CR-equivalent to $\6T$ is the universal covering
of a domain in $\6R$ on which a suitable subgroup of the
CR-automorphism group $\Aut(\6R)$ acts transitively. What concerns
some natural group actions, we observe that the nonclosed orbits in
the complex manifold $\SL(2,\CC)$ under the action of the group
$\SL(2,\RR)\times\SL(2,\RR)$ given by $z\mapsto gzh^{-1}$ or the
nonclosed orbits in the same complex manifold, but under a different
action, namely of $\SL(2,\CC)\cong \Spin(1,3)$ acting by $z\mapsto
gz\overline g'$, are all CR-equivalent to certain homogeneous
domains in $\6R$.

We close by indicating in Section \ruf{final} how some results from
Section \ruf{Global} can be generalized to higher dimensions.
}

\KAP{Preliminaries}{Preliminaries} 

In this paper we use essentially the same conventions and notation as
in \LIT{KAZT}{section 2}: Let $Z$ be a complex manifold $Z$ and
$\pi:TZ\to Z$ its tangent bundle. Then $TZ$ also has the structure of
a complex manifold with $\pi$ being a holomorphic submersion. In
particular, every tangent space $T_{a}Z$, $a\in Z$, is a complex
vector space with $\dim_{\CC}T_aZ=n$ if $Z$ has complex dimension
$n$ at $a$. By $\hol(Z)$ we denote the complex Lie algebra of all
holomorphic vector fields on $Z$, that is, of all holomorphic sections
$\xi:Z\to TM$ in the tangent bundle over $Z$. For every $a\in Z$ we
denote the corresponding tangent vector in $T_{a}Z$ with $\xi_{a}$
(and not with $\xi(a)\,$).

A (locally closed) real submanifold $M\subset Z$ is called a {\sl
CR-manifold} if the complex dimension of the {\sl holomorphic tangent
space} $H_{a}M:=T_{a}M\cap iT_{a}M$ is a locally constant function of
$a\in M$ (here every tangent space $T_{a}M$ is considered as an
$\RR$-linear subspace of $T_{a}Z$). A smooth function $f:M\to\CC$ is
called CR if for every $a\in M$ the restriction of its differential to
$H_{a}M$ is complex linear. More generally, a smooth mapping
$\phi:M\to M'$ between CR-manifolds is CR, if the differential
$d_{a}\phi:T_{a}M\to T_{\phi(a)}M'$ maps $H_{a}M$ to $H_{\phi(a)}M'$
and is complex linear thereon. The notion of CR-function and
CR-mapping can also be generalized to the nonsmooth case; then the
conditions have to hold in the distribution sense, see \Lit{BERO} as a
general reference for CR-manifolds.

In the following we only consider connected real-analytic
CR-submanifolds $M$ of $Z$. With $\hol(M)$ we denote the real Lie
algebra of all vector fields $\eta:M\to TM\subset TZ$ on $M$ with the
following property: To every $a\in M$ there is an open neighbourhood
$U$ of $a$ with respect to $Z$ and a holomorphic vector field
$\xi\in\hol(U)$ with $\xi_{z}=\eta_{z}$ for all $z\in U\cap
M$. Clearly, in case the CR-submanifold $M$ is {\sl generic} in $Z$
(that is, satisfies $T_{a}Z=T_{a}M+iT_{a}M$ for every $a\in M$), the
local holomorphic extension $\xi$ above can always be chosen in such a
way that $U$ is an open neighbourhood of all of $M$ in $Z$. The
elements of $\hol(M)$ are also called {\sl infinitesimal
CR-transformations} on $M$, the reason being that the corresponding
local flow on $M$ consists of real-analytic CR-transformations. The
vector field $\eta\in\hol(M)$ is called {\sl complete on } $M$ if the
local flow actually consists of a one-parameter family $(g_{t})$,
$t\in\RR$, of global transformations. Then we write
$\exp(\eta):=g_{1}$ and get a mapping $\exp:\aut(M)\to\Aut(M)$, where
$\aut(M)\subset\hol(M)$ is the subset of all complete vector fields
and $\Aut(M)$ is the group of all real-analytic CR-diffeomorphisms of
$M$. In general, $\aut(M)\subset\hol(M)$ is not closed under addition
nor under taking brackets. But in case there is a Lie subalgebra
$\7g\subset\hol(M)$ of finite dimension containing $\aut(M)$, then a
result of \Lit{PALA} implies: $\aut(M)$ itself is a Lie algebra of
finite dimension and $\Aut(M)$ has a unique Lie group structure such
that $\exp:\aut(M)\to\Aut(M)$ is locally bianalytic in a neighbourhood
of $0\in\aut(M)$.

For every $a\in M$ denote by $\hol(M,a)$ the space of all germs at $a$
of infinitesimal CR-trans\-for\-mations defined in arbitrary open
neighbourhoods of $a$. Then $\hol(M,a)$ is a real Lie algebra, and
$\aut(M,a):=\{\xi\in\hol(M,a):\xi_{a}=0\}$ is a Lie subalgebra of
finite codimension. We call the CR-manifold $M$ {\sl locally
homogeneous} if the evaluation map $\hol(M,a)\to T_{a}M$,
$\xi\mapsto\xi_{a}$, is surjective for every $a\in M$. Local
homogeneity implies that to every pair $a,b$ of points in $M$ there
exist open neighbourhoods $U,V$ of $a,b$ in $M$ together with a
real-analytic CR-diffeomorphism $U\to V$. We call $M$ {\sl
homogeneous} if there exists a connected Lie group $G$ together with a
group homomorphism $\Phi:G\to\Aut(M)$ such that the mapping $G\times
M\to M$, $\,(g,x)\mapsto \Phi(g)x$, is real-analytic and $G$ acts
transitively on $M$ via $\Phi$. Clearly, `homogeneous' implies
`locally homogeneous'.

The CR-manifold $M\subset Z$ is called {\sl holomorphically
nondegenerate} if for every domain $U\subset Z$ and every
$\xi\in\hol(U)$ with $\xi_{x},i\[1\xi_{x}\in T_{x}M$ for all $x\in
U\cap M$ necessarily $\xi_{x}=0$ holds for all $x\in U\cap M$. In case
$M$ is a real-analytic hypersurface in $Z$, holomorphic nondegeneracy
is equivalent to $\dim\hol(M,a)<\infty$ for all $a\in M$, see
e.g. \Lit{BERO} \p367 for this and related results. In this note
(except in the final section \ruf{final}) the complex manifold $Z$
always has dimension $3$ (either $\CC^{3}$ or a nonsingular quadric in
$\PP_{4}(\CC)\,$) and $M$ is a real hypersurface that is
2-nondegenerate, a property that implies `holomorphically
nondegenerate'.

\Ues{Convention for notating vector fields} In this paper we do not
need the complexified tangent bundle $TM\otimes_{\RR}\CC$ of $M$. All
vector fields occurring here correspond to `real vector fields'
elsewhere. In particular, if $E$ is a complex vector space of finite
dimension and $U\subset E$ is an open subset then the vector fields
$\xi\in\hol(U)$ correspond to holomorphic mappings $f:U\to E$, and the
correspondence is given in terms of the canonical trivialization
$TU\cong U\Times E$ by identifying the mapping $f$ with the vector
field $\xi=(\id_{U},f)$. To have a short notation we also write
$$\xi=f(z)\ddz{}\,.$$ As an example, if $E$ is the space of all
complex $n\Times m$-matrices and $c$ is an $m\times n$-matrix, then
$zcz\ddz{}$ denotes the quadratic vector field on $E$ corresponding to
the holomorphic mapping $E\to E$, $z\mapsto zcz$. As soon as the
vector field $\xi=f(z)\ddz{}$ is considered as differential operator,
special caution is necessary: $\xi$ applied to the smooth function $h$
on $U$ is
$\xi h=f(z)\ddz{}h+\overline f(z)\ddo h$. We therefore stress
again that we write
$$\xi=f(z)\ddz{}\Steil{ instead of }\xi=f(z)\ddz{}+\overline
f(z)\ddo{}\hbox{~~elsewhere}\,,\Leqno{SF}$$ and this convention will
be in effect allover the paper. In particular, in case $E=\CC^{n}$
with coordinates $z=(z_{1},\dots,z_{n})$ we write
$$\xi=f_{1}(z)\ddz{1}+f_{2}(z)\ddz{2}+\dots+f_{n}(z)\ddz{n}\,,$$ where
$f=(f_{1},\dots,f_{n}):U\to\CC^{n}$ is holomorphic.

\KAP{tube}{The tube over the light cone} 

\beginabs In this section we introduce the light cone tube $\6T$, a
real hypersurface of $\CC^{3}$ that is everywhere Levi degenerate but
has finite dimensional Lie algebra $\hol(\6T,a)$ at every point
$a\in\6T$. From \Lit{KAZT} we recall the explicit form of this Lie
algebra and present a root decomposition in terms of vector fields.
This will enable us in the next section to show that a certain local
equation defines a CR-manifold locally equivalent to $\6T$.\endabs

A convenient model for $3$-dimensional space time is the linear
subspace $V\subset\RR^{2\Times2}$ of all real symmetric
$2\times2$-matrices with the (normalized) trace as time coordinate.
There
$$\Omega:=\{v\in V:v\hbox{ positive definite}\}$$ is the {\sl
future cone} and its smooth boundary part is the {\sl future light
cone}
$$\6C\colon=\big\{v\in V:\det(v)=0\,,\tr(v)>0\big\}=\left\{
\Matrix{t+x_{1}}{x_{2}}{x_{2}}{t-x_{1}} \in
V:t^{2}=x^{2}_{1}+x^{2}_{2}\,,\;t>0\right\}\,.$$ It is obvious that there
exists a $2$-dimensional group of linear transformations on $\RR^{3}$
acting simply transitive on $\6C$. 

\smallskip The main object of our interest is the tube
$$\6T\colon=\6C\oplus iV=\{z\in V\oplus iV:\det(z+\overline
z)=0\;,\Re\tr(z)>0\}\Leqno{DF}$$ over $\6C$, where we identify the
complexification $V\oplus iV$ in the obvious way with the space $E$ of
all symmetric complex $2\Times2$-matrices. $\6T$ is the smooth
boundary part of the {\sl tube domain} $\6H\colon=\Omega\oplus iV$ in
$E$ ({\sl Siegel's upper half plane} up to the factor $i$) and is a
locally closed real-analytic hypersurface of $E$ with everywhere
degenerate Levi form. Actually it is well known that $ \6T$ is
everywhere $2$-nondegenerate as CR-manifold.

\Ues{The CR-automorphisms} The 7-dimensional Lie
group of affine transformations on $E$
$$\{z\mapsto gzg'+iv:g \in \GL(2,\RR), v\in V\}\Leqno{BC}$$ acts
transitively on $\6T$ and $\6H$, where $g'$ denotes the transpose
of $g$. It is known \Lit{KAZT} that $\Aut(\6T)$ is the group of all
transformations \Ruf{BC} while $\Aut(\6H)$ is the $10$-dimensional
group of all biholomorphic transformations
$$z\;\mapsto\;(az-ib)(icz+d)^{-1}\,,\Leqno{WC}$$ where $\Matrix abcd$
is in the real symplectic subgroup $\Sp(2,\RR)\subset\SL(4,\RR)$ with
$a,b,c,d\in\RR^{2\Times2}$, see \Lit{KLIN} \p351. Differentiating the
action of $\Sp(2,\RR)$ gives
$$\7g:=\aut(\6H)=\big\{(b+cz+zc'+zdz)\ddz{}:b,d\in
iV,\,c\in\RR^{2\Times2}\big\}\cong\7{sp}(2,\RR)\cong\7{so}(2,3)\Leqno{FS}$$
with convention \Ruf{SF} in effect. All vector fields in $\aut(\6H)$
are polynomial of degree $\le2$ on $E$, in particular,
$$\aut(\6H)\subset\hol(\6T)\subset\hol(\6T,a)$$ in a canonical way for
all $a\in\6T$. One of the main results of \Lit{KAZT} states (even for
higher dimensional examples of this type) that all three Lie algebras
above coincide, see Proposition 4.3 in \Lit{KAZT}.

\smallskip The vector fields in $\7g$ corresponding to
$c=\Matrix1000,\Matrix0001$ in \Ruf{FS} are
$$\zeta_{1}:=2z_{0}\ddz0+z_{1}\ddz1\Steil{and}\zeta_{2}:=z_{1}\ddz1+2
z_{2}\ddz2\,,\Leqno{FX}$$ when expressed in the coordinates
$$(z_{0},z_{1},z_{2})\longmapsto\pmatrix{z_{0}&z_1\cr z_{1}& z_{2}\cr}$$ on
$E$. These give the following decomposition 
$$\7g=\bigoplus_{\nu\in\ZZ^{2}}\7g^{\nu}\Steil{with}[\7g^{\mu},\7g^{\nu}]
\subset\7g^{\mu+\nu}\;,\Leqno{OR}$$ where
$\7g^{\nu}:=\{\xi\in\7g:[\zeta_{j},\xi]=\nu_{j}\[1\xi\steil{for}j=1,2\}$
for $\nu=(\nu_{1},\nu_{2})$ and, in particular,
$\7g^0=\RR\zeta_{1}\oplus\RR\zeta_{2}$. For every root (i.e. $\nu\ne0$
and $\7g^{\nu}\ne0$) the corresponding root space $\7g^{\nu}$ has real
dimension $1$, and the set of all roots is visualized by the eight
vectors in {\tt Figure $\!$1}, the root system of the complex simple
Lie algebra $\7{so}(5,\CC)\cong\7{sp}(2,\CC)$. \midinsert
$$\epsfxsize=40 true mm
 \hbox{\rlap{\hskip95pt\raise33pt\hbox{$\nu^{}_{1}$}}
 \rlap{\hskip52pt\raise73pt\hbox{$\nu^{}_{2}$}}
 \epsffile{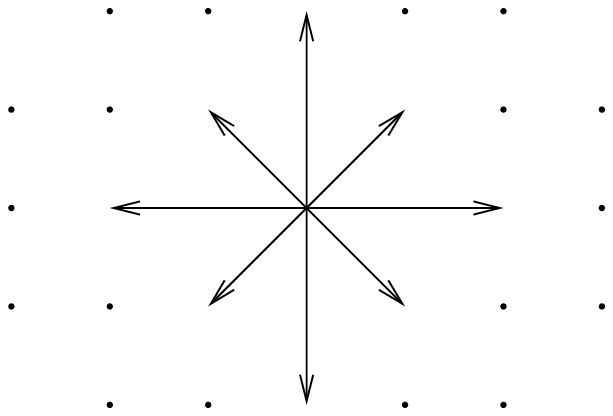}}\hskip2cm \epsfxsize=35 true mm
 \hbox{\epsffile{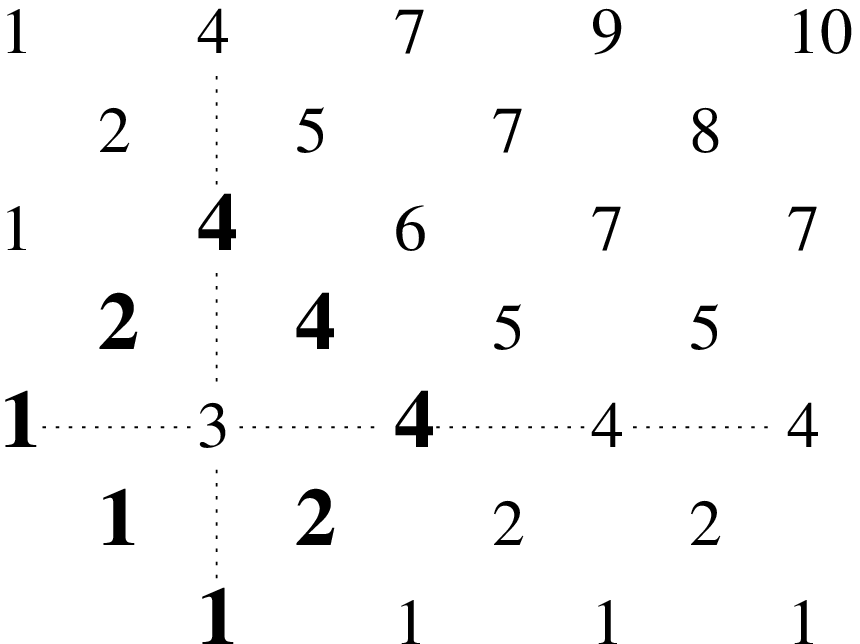}}$$ \hskip4.5 true cm{\eightt Figure
 1\hskip4.6truecm Figure 2}
\endinsert

\noindent An explicit choice of root vectors $\xi^{\nu}\in\7g^{\nu}$ 
is as follows:

{\Klein\bigskip\advance\baselineskip5pt
\line{\rlap{\Write{DR}{\tenrm(\PrN)}}\hfill$\xi^{0,2}=
iz_{1}^{2}\ddz0+iz_{1} z_{2}\ddz1+i z_{2}^{2}\ddz2$\hfill}
\centerline{\qquad\qquad$\xi^{-1,1}=2z_1\ddz0+
z_{2}\ddz1\hfill\xi^{1,1}=2iz_{0}z_{1}\ddz0+i(z_{1}^{2}+z_{0}
z_{2})\ddz1+2iz_{1}z_{2}\ddz2$\qquad}
\centerline{$\xi^{-2,0}=i\ddz0\hfill\xi^{2,0}=
iz_{0}^{2}\ddz0+iz_{0}z_1\ddz1+iz^{2}_{1}\ddz2$}
\centerline{\qquad\qquad$\xi^{-1,-1}=i\ddz1\qquad\qquad\qquad\qquad
\xi^{1,-1}=-z_{0}\ddz1-2z_1\ddz2$\qquad\qquad}
\centerline{$\xi^{0,-2}=i\ddz2\,.$}}

\Ues{Remarks} The choice of the vector fields
\Ruf{FX} gives the Cartan subalgebra
$\7g^{0}=\RR\zeta_{1}\oplus\RR\zeta_{2}$ of $\7g$. Unlike to the
complex situation there may exist several nonconjugate Cartan
subalgebras in a real semi-simple Lie algebra. In our particular
situation $\7g\cong\7{so}(2,3)$ there are precisely four conjugacy
classes of Cartan subalgebras.

Note that all vector fields in $\7g^{0}$ vanish simultaneously in a
unique point of the nonsmooth boundary part $iV$ of $\6H$, namely in
the origin $0\in E$. In Section \ruf{local} we will consider another
Cartan subalgebra of $\7g$ that vanishes in a point of the smooth
boundary part $\6T$ of $\6H$, and in Section \ruf{bounded} even a
third one is considered that vanishes in a point of $\6H$,
compare \Ruf{TR}.

It is not a coincidence that the Cartan subalgebras above vanish at
some point of $E$. To put it in a broader perspective, $E$ admits a
natural embedding in a compact complex flag manifold $Z$ together with
a global action of $\Sp(2,\RR)$ on $Z$ that extends the action
\Ruf{WC}. By Borel's fixed point theorem every torus $H:=\exp(\7h)$,
$\7h\subset\7g$ a Cartan subalgebra, has a fixed point in $Z$.

The only nonclosed orbits in $E$ of the group $$H:=\{z\mapsto
gzg'+iv:g\in\SL(2,\RR), v\in V\}$$ are $\6T$ and $-\6T$. Furthermore,
the real function $\det(z+\overline z)$ is constant on $H$-orbits and
its sign determines the signature of the Levi form for the orbits.

\KAP{local}{A local realization}

\beginabs In this section we introduce a family $(\6M_{t})_{t\in\RR}$ of local 
CR-submanifolds in $\CC^{3}$ that may be considered as local
deformations of the tube $\6T$, are pairwise CR-inequivalent and have
semi-simple Lie algebras $\hol(\6M_{t},0)$, $t\ne0$, of dimension
$\ge6$. In particular, with $\6M_{1}$ we get a normalized equation for
$\6T$ in the sense of \Lit{EBEN}.\endabs

In the following we consider in $\CC^{3}$ with coordinates
$(w,z_{1},z_{2})$ hypersurfaces $M$ of the following type. $M$ is
given near $0\in\CC^{3}$ by an equation
$$ w+\overline w\;=\;2z_{1}\overline z_{1}+\big(z_{1}^2\overline
z_{2}+\overline z_{1}^2z_{2}\big)+q\big(\Im(w),z_{1},\overline
z_1,z_{2},\overline z_2\big)\;,\Leqno{AE}$$ where $q$ is a convergent
real power series whose nonzero monomial terms either have degree
$\ge4$ or have degree $3$ and then contain $\Im(w)$. Then \Ruf{AE} is
just the partial normal form (A,ii,2) in \Lit{EBFT}. In particular, $M$ is
$2$-nondegenerate at $0$ and there the Levi kernel is the
$z_{2}$-coordinate axis. Further, $\dim \hol(M,0)<\infty$ holds.

Consider on $\CC^3$ the group of all linear transformations
$(w,z_{1},z_{2})\mapsto(s^{2}w,stz_{1},t^{2}z_{2})$ with $s\in
e^{\RR}$ and $t\in e^{i\RR}$. Clearly, this group is the exponential
of the Lie algebra spanned in $\aut(\CC^{3})$ by the two (real)
holomorphic vector fields $\zeta_1$ and $i\zeta_2$ where, using
convention \Ruf{SF}, $\zeta_1$ and $\zeta_2$ are given in the local
coordinates by
$$\zeta_{1}:=2w\ddw+z_{1}\!\ddz1\Steil{and}\zeta_{2}:=
z_{1}\!\ddz1+2z_{2}\ddz2\,.\Leqno{AF}$$ The coordinate expression of
these vector fields coincides with those in \Ruf{FX} (if the variable
$w$ is renamed $z_{0}$). The same will occur in Section \ruf{bounded}.

From now on we assume for the rest of the section that
$\7g:=\hol(M,0)$ contains the vector fields $\zeta_1$ and $i\zeta_2$.
For the sake of clarity let us emphasize that under this assumption
the vector field $\zeta_2$ cannot be contained in $\7g$ since $M$ is
holomorphically nondegenerate.

\medskip Our first goal is to show that under the assumption
$\zeta_{1},i\[1\zeta_{2}\in\7g$ there is a unique equation of the form
\Ruf{AE} such that $M$ is locally homogeneous at $0\in M$. It will
turn out that in this case $M$ is locally CR-equivalent to the tube
$\6T$ over the future light cone in $\CC^{3}$ as considered in Section
\ruf{Preliminaries}.

We need some preparation: Let $\7P$ be the complex Lie algebra of all
polynomial holomorphic vector fields on $\CC^{3}$ with coordinates
$w,z_{1},z_{2}$. The natural adjoint action of $\zeta_1$ and $i\zeta_2$
induces a grading of $\7P$: Consider the lattice
$$\Lambda:=\{n+im\in\ZZ+i\ZZ:n+m\in2\ZZ\}\Leqno{AK}$$ in $\CC$ and
denote by $\7P^{\lambda}$ the $\lambda$-eigenspace of
$\ad(\zeta_{1}+i\[1\zeta_{2})$ in $\7P$. Note that every monomial
vector field is contained in some $\7P^{\lambda}$: For all
$m,n,l\in\NN$ and $k=0,1,2$ the vector field
$w^{m}\]2z_{1}^{n}z_{2}^{l}\ddz k$ (with $\ddz 0:=\ddw$ to simplify
notation) is in $\7P^{\lambda}$ for $\lambda=(2m+n+k-2)+i(n+2l-k)$. In
particular, every $\7P^{\lambda}$ has finite dimension and
$$\7P=\bigoplus_{\lambda\in\Lambda}\7P^{\lambda}\,,\qquad
[\7P^{\lambda},\7P^{\mu}]\subset\7P^{\lambda+\mu}\,.\Leqno{AT}$$ A
small check shows that $\7P^{\lambda}=0$ if
$\min(\Re\lambda,\Im\lambda,\Re\lambda+\Im\lambda)<-2$. For instance,
for $\lambda\in \{-2,\,-1-i,\,-2i\}$ the spaces
$\7P^{\lambda}=\CC\xi^{\lambda}$ are 1-dimensional with generators
$$\xi^{-2}:=i\ddw,\quad\xi^{-1-i}:=
\ddz1,\quad\xi^{-2i}:=\ddz2\,.\Leqno{AH}$$ In {\tt Figure $\!\!$2$\;$}
the nonzero complex dimensions of $\7P^{\lambda}$ are listed for all
$\lambda=m+in$ with $m\le6$ and $n\le4$.

We now explain how $\7g$ is related to the decomposition
$\7P=\bigoplus \7P^\lambda$. A priori, the finite dimensional Lie
algebra $\7g$ is contained in $\hol(U)$ for some open neighbourhood $U$ of
$0\in\CC^{3}$. The sum $\7l:=\7g+i\7g$ in $\hol(U)$ actually is a
direct sum of real subspaces. Since the complex Lie algebra $\7l$
contains the Euler vector field $(\zeta_{1}+\zeta_{2})/2$ necessarily
$\7l$ is contained in $\7P$. Both $\7l$ and $\7g$ are invariant under
$\ad(\zeta_{1}+i\[1\zeta_{2})$. This gives immediately the
decomposition
$$ \7l=\bigoplus_{\lambda\in\Lambda}\7l^{\lambda},\qquad
\7l^{\lambda}:=\7l\cap\7P^{\lambda}\,.$$ Let $\xi\mapsto\overline\xi$
be the conjugation of $\7l$ with respect to the real form $\7g$. Then
$$\overline{\7l^\lambda}=\7l^{\overline\lambda}\Leqno{CI}$$ for all
$\lambda$ and consequently $\7g\cap\7l^\lambda=0$ if
$\lambda\not\in\RR$. Furthermore, $\7g$ admits the generalized
eigenspace decomposition
$$\7g=\bigoplus_{\hbox{\sevenrm
Im}\lambda\ge0}\7g^{[\lambda]}\,,\qquad\7g^{[\lambda]}:=
\7g\cap\big(\7l^{\lambda}+ \7l^{\overline\lambda}\big)\,.\Leqno{AG}$$

\medskip Next we claim that
$$\7g^{[0]}=\RR\zeta_{1}\oplus\RR i\[1\zeta_{2}\Leqno{AY}$$ holds.
Indeed, assume to the contrary that there exists $\xi\in\7g^{[0]}$
with $\xi\notin(\RR\zeta_{1}\!+\RR i\[1\zeta_{2})$. Then
$\xi\in\7P^{0}$ and, after subtracting a suitable linear combination
of $\zeta_{1}$, $i\zeta_{2}$, we may assume $\xi=\alpha w\ddw+\beta
z_{2}\ddz2$ for certain $\alpha,\beta\in\CC$. Applying $\xi$ to the
defining equation \Ruf{AE} yields on $M$ for $r:=\Re(\alpha)$ and
$s:=\Im(\alpha)$ the identity $$ r\[2(2z_{1}\overline
z_{1}+z_{1}^2\overline z_{2}+\overline
z_{1}^2z_{2})-2s\[2\Im(w)=\beta\,\overline z^{2}_{1}z_{2}+\overline
\beta\,z^{2}_{1}\overline z_{2}+\cdots$$ up to terms of degree $\ge4$
or of degree $3$ containing $\Im(w)$ (the convention \Ruf{SF} has to
be observed). This forces $\alpha=\beta=0$ and, in particular,
$\xi\in(\RR\zeta_{1}\!+\RR i\[1\zeta_{2})$.\qed

\vfill\eject

Our first result now states 

\Proposition{AR} Suppose that the CR-manifold $M$ is given locally by
the equation \Ruf{AE} and that $\7g=\hol(M,0)$ contains the vector
fields $\zeta_{1},i\[1\zeta_{2}$. Then the following conditions are
equivalent. \0 $M$ is locally homogeneous at $0$. \1 $\7g$ contains
the vector fields\hfill\break $\eta:=2z_{1}\ddw+(1-z_{2})\ddz1$ and
$\chi:=z_{1}^{2}\ddw -z_{1}z_{2}\ddz1+(1-z^{2}_{2})\ddz2\,.$ \1 The
term $q$ in the defining equation \Ruf{AE} is\nline
\phantom{XXXXXXXXXXXXX} $q=(2z_{1}\overline z_{1}+z_{1}^2\overline
z_{2}+\overline z_{1}^2z_{2})\cdot\sum_{j=1}^\infty(z_2\overline
z_2)^{j}$,\nline that is, the defining equation reads $w+\overline
w=(2z_{1}\overline z_{1}+z_{1}^2\overline z_{2}+\overline
z_{1}^2z_{2})(1-z_{2}\overline z_{2})^{-1}$.

\Proof \To12 Suppose that $M$ is locally homogeneous at $0$. Then
$\{\xi_{0}:\xi\in\7l\}=T_0M+i\[1T_0M\cong\CC^3$. Since $\7l^{\lambda}$
vanishes at $0$ for every $\lambda\notin\{-2,\,-1-i,\,-2i\}$ we get
$\dim\7l^{\lambda}=1$ for all $\lambda\in\{-2,\,-1\pm i,\,\pm2i\}$
where we have used \Ruf{CI}. Since every $\7g^{[\lambda]}$ is
invariant under $\ad(i\zeta_{2})$ we have
$\{\xi_{0}:\xi\in\7g^{[\lambda]}\}= \{\xi_{0}:\xi\in\7l^\lambda\}$ if
$\Im(\lambda)<0$. Therefore there exists a unique $\xi\in\7P^{-1+i}$
with $\eta:=\xi^{-1-i}+\xi\in\7g$, that is, $\eta=\ddz1+\alpha
z_{1}\ddw+\beta z_{2}\ddz1\in\7g$ for certain $\alpha,\beta\in\CC$.
Applying $\eta$ to the equation \Ruf{AE} and comparing homogeneous
terms immediately gives $\alpha=2$ and $\beta=-1$. In the same way we
get a vector field $\chi:=\ddz2+z_{1}^{2}\ddw -\delta
z_{1}z_{2}\ddz1-\epsilon z^{2}_{2}\ddz2\in\7g^{[2i]}$ for suitable
$\delta,\epsilon\in\CC$. But then
$[\7g^{[-1+i]},\7g^{[2i]}]\subset\7g^{[-1+i]}$ and
$$[\eta,\chi]=\big(2z_{1}+2(\delta-1)z_{1}z_{2}\big)\ddw+\big(1-\delta
z_{2}+(\delta-\epsilon)z_{2}^{2}\big)\ddz1\;\in\;\7g^{[-1+i]}\,,$$
imply $[\eta,\chi]=\eta$, i.e. $\delta=\epsilon=1$. \nline \To21 In
case $\eta\in\7g^{[-1+i]}$, $\chi\in\7g^{[2i]}$ the linear subspace
$\7g^{[-1+i]}\oplus\7g^{[2i]}\subset\7g$ spans the holomorphic tangent
space $H_{0}(M)\subset T_{0}M$. The vector field
$\big[\eta,[i\zeta_{2},\eta]\big]\in\7g^{[-2]}$ does not vanish and
thus spans $T_{0}M/H_{0}M$, that is, $\7g$ spans the full tangent
space $T_{0}M$. \nline \To32 This is easily checked. \nline \To23
The $\RR$-linear span of
$\xi^{-2},\,\eta,\,\chi,\,[i\zeta_{2},\eta],\,
[i\zeta_{2},\chi],\,i\zeta_{2}$ is a Lie subalgebra $\7a\subset\7g$ of
dimension $6$ which also spans the full tangent space of $M$ at
$0$. Therefore $M$ is the the local integral manifold of $\7a$ near
$0\in\CC^{3}$. Denote by $\tilde M$ the hypersurface defined near
$0\in\CC^{3}$ by the equation \Ruf{AE}, where $q$ is replaced by
$\tilde q=(2z_{1}\overline z_{1}+z_{1}^2\overline z_{2}+\overline
z_{1}^2z_{2})\cdot\sum_{j=1}^\infty(z_2\overline z_2)^{j}$. By a
routine check it is verified that also $\7a\subset\hol(\tilde M,0)$
holds, that is, also $\tilde M$ is the local integral manifold of
$\7a$ near $0\in\CC^{3}$. This implies that the germs of $M$ and
$\tilde M$ at $0$ coincide and hence that $q=\tilde q$ as power
series.\qed

\Ues{A family of CR-manifolds} In the following we fix $t\in\RR$ and
consider the CR-hypersurface $\6M_{t}$ given in $\CC^{3}$ near the
origin by the equation
$$w+\overline w=(2z_{1}\overline z_{1}+z^{2}_{1}\overline
z_{2}+\overline z^{2}_{1}z_{2})(1-tz_{2}\overline
z_{2})^{-1}\,.\Leqno{AJ}$$ Clearly, this is a special case of
\Ruf{AE}, and it is easily checked that the vector fields
$\zeta_1,i\zeta_2$ are contained in
$\7g:=\7g_{t}:=\hol(\6M_{t},0)$. The manifold $\6M_{t}$ for $t=0$ has
already been studied in \Lit{ERSH} and for $t=1$ in \Lit{GAME}, where
also the vector fields $\eta$, $\chi$ from \Ruf{AR} occur.

Our next goal is to describe how the Lie algebra $\7g$ depends on the
parameter $t$. Put $\Psi:=\{\pm2,\,\pm2i\}$ and define for each
$\lambda\in\Psi\subset\Lambda$ the vector fields
$\xi^{\lambda}\in\7P^{\lambda}$ as follows:
$$\xi^{2i}:=
z_{1}^{2}\ddw-tz_{1}z_{2}\ddz1-tz^{2}_{2}\ddz2\ , \qquad
\xi^{2}:=itw^{2}\ddw+itwz_{1}\ddz1 -iz_{1}^{2}\ddz2\;\Leqno{AZ}$$
and $\xi^{-2i},\xi^{-2}$ as in \Ruf{AH}. 
Furthermore put
$$\eta^{\lambda}:=\cases{\xi^{\lambda}+\xi^{\overline\lambda}&
$\Im(\lambda)\ge0$\cr i\xi^{\lambda}-i\[1\xi^{\overline\lambda}&
otherwise\cr}\Leqno{AI}$$ and notice that
$$\eta^{\lambda}\in\7g,\qquad[\zeta_{1},\eta^{\lambda}]=
\Re(\lambda)\,\eta^{\lambda}\Steil{and}
[i\[1\zeta_{2},\eta^{\lambda}]=\Im(\lambda)\,\eta^{\overline\lambda}$$
hold for every $\lambda\in\Psi$. By applying vector fields to equation
\Ruf{AJ} it is not difficult to see that the Lie algebra
$$\eqalign{\7h:&=\;(\RR\eta^{-2}\oplus\RR\zeta_{1}\oplus\RR\eta^{2})\oplus
(\RR\eta^{-2i}\oplus\RR i\[1\zeta_{2}\oplus\RR\eta^{2i})\cr &\cong\;\cases{
\7{sl}(2,\RR)\times\7{sl}(2,\RR)&$t>0$\cr
(\7{so}(1,1)\ltimes\RR^{2})\times(\7{so}(2)\ltimes\RR^{2})&$t=0$\cr
\7{sl}(2,\RR)\times\7{su}(2)&$t<0$\cr }\cr}\Leqno{AW}$$ is contained
in $\7g$. In case $t\ne1$ the CR-manifold $\6M_{t}$ cannot be locally
homogeneous at $0$ since then the vector field $\chi$ from \Ruf{AR}
is not contained in $\7g^{[2i]}=\RR\eta^{2i}\oplus\RR\eta^{-2i}$,
compare Proposition \ruf{AR}. In case $t=0$ the right hand side of
\Ruf{AJ} reduces to a cubic polynomial and, with some computation, it
can be seen that $\7g=\7h$ in this case, compare also case B2 in
\LIT{ERSH}{\p94}. In all other cases $\7h$ is semi-simple. We shall use
this in the following.

\Ues{Explicit determination of
$\,{\bfrak{h\]2o\]1l}\3{\9(}\9{\6M}_{\3t},{\bf 0}\9)}\,$} ~For fixed
$t\in\RR$ and $\7l=\7g\oplus i\7g$ with $\7g=\hol(\6M_{t},0)$ put
$$\Phi:=\{\lambda\in\Lambda:\lambda\ne0\steil{and}\7l^{\lambda}\ne0\}\,.$$
Then $\Psi\subset\Phi$ and $\Phi$ is invariant under
$\lambda\mapsto\overline\lambda$. For every $k\in\ZZ$ denote by
$d_{k}$ the dimension of the eigenspace
$\{\xi\in\7l:[\zeta_{1},\xi]=k\xi\}=\bigoplus_{\hbox{\sevenrm
Re}\lambda=k}\7l^{\lambda}$. In case $t\ne0$ the representation theory
of $\7{sl}(2,\CC)\cong\CC\eta^{-2}\oplus\CC\zeta_{1}\oplus\CC\eta^{2}$
implies $d_{k}=d_{-k}$ for all $k$, and thus $d_{1}\in\{0,2\}$,
$d_{2}=1$ and $d_{k}=0$ for $k>2$. As a consequence, $\Phi$ is also
invariant under $\lambda\mapsto-\overline\lambda$ and $\7l^{\lambda}$
has dimension $1$ for every $\lambda\in\Phi$. Because of
$(-1-i)\notin\Phi$ in case $t\ne1$ we therefore have $\Phi=\Psi$ in
this situation. This proves:

\Proposition{} In case $t\ne1$ the Lie algebra $\hol(\6M_{t},0)$
coincides with $\7h$ and, in particular, has dimension $6$.\Formend

\medskip It remains to consider the case $t=1$. Applying $\eta$ to
equation \Ruf{AJ} implies $\eta\in\7g=\hol(\6M_{1},0)$. In particular
$(-1-i)\in\Phi$ and hence $\Phi:=\{\pm2,\,\pm1\pm i,\,\pm2i\}$, that
is, $\Phi$ can also be visualized by {\tt Figure 1}, while
$\Psi\subset\Phi$ corresponds to the subset of long arrows. Define for
every $\lambda\in\Phi$ the vectors $\xi^{\lambda}\in\7P^{\lambda}$ by
$\Ruf{AH}$ and
$$\eqalign{\xi^{-1+i}:&=2z_{1}\ddw-z_{2}\ddz1\cr\xi^{1-i}:&=
2z_{1}\ddz2-w\ddz1\cr} \quad\eqalign{\xi^{2i}:&=
z_{1}^{2}\ddw-z_{1}z_{2}\ddz1-z^{2}_{2}\ddz2\cr
\xi^{1+i}:&=2wz_{1}\ddw+ (z^{2}_{1}-wz_{2})\ddz1+2z_{1}z_{2}\ddz2
\cr\xi^{2}:&=iw^{2}\ddw+iwz_{1}\ddz1
-iz_{1}^{2}\ddz2\;.\cr}\Leqno{AL}$$ Notice that after replacing
$z_{2}$ by $-z_{2}$ all vector fields in \Ruf{AL} and \Ruf{AH} become
complex multiples of those in \Ruf{DR}.

Now define $\eta^{\lambda}$ by \Ruf{AI} for every
$\lambda\in\Phi$. Differentiating the defining equation for $M_{1}$
along $\eta^{\lambda}$ gives $\eta^{\lambda}\in\7g$ and hence
$\7l^{\lambda}=\CC\xi^{\lambda}$ for all $\lambda\in\Phi$. For all
$\lambda,\mu\in\Phi$ with $\lambda+\mu\ne0$ the identity
$[\7l^{\lambda},\7l^{\mu}]=\7l^{\lambda+\mu}$ can be verified, that
is, $\7l$ is isomorphic to $\7{so}(5,\CC)$. As real form of $\7l$
therefore $\7g$ is isomorphic to $\7{so}(5)$, $\7{so}(1,4)$ or
$\7{so}(2,3)$. The isotropy subalgebra $\7l_{0}$ of $\7l$ at the
origin of $\CC^{3}$ is maximal parabolic. Therefore, $\6M_{1}$ occurs
as piece of an $\R$-orbit in a 3-dimensional quadric
$Z\subset\PP_{4}(\CC)$ or in $\PP_{3}(\CC)$, where $\R$ is a real form
of $\SO(5,\CC)\cong\Sp(2,\CC)$. Since the only $2$-nondegenerate
hypersurface orbits occur in case of $\SO(2,3)$ acting on $Z$, the Lie
algebra $\7g$ can only be isomorphic to $\7{so}(2,3)$. Since the tube
$\6T$ over the future light cone actually can be realized as an open
piece of an $\SO(2,3)$-orbit in $Z$, we get (a direkt proof is given
by formula \Ruf{OC}):

\Proposition{IL} The hypersurface $\6M_{1}$ defined near $0\in\CC^{3}$ by
the equation $$w+\overline w=(2z_{1}\overline z_{1}+z^{2}_{1}\overline
z_{2}+\overline z^{2}_{1}z_{2})(1-z_{2}\overline
z_{2})^{-1}\Leqno{AO}$$ is locally CR-isomorphic to the tube $\6T$ over
the future light cone. \Formend

The local CR-equivalence between $\6M_{1}$ and $\6T$ can be made
explicit. For this, a transformation in $\SO(5,\CC)$ has to be found
that the corresponding $\SO(2,3)$-orbits in the quadric $Z$ maps
to each other. Let $U:=\{(w,z_{1},z_{2})\in\CC^{3}:|z_{2}|<1\}$ and
consider $\6M_{1}$ as closed CR-submanifold of $U$ via the equation
\Ruf{AO}. Then
$$(w,z_{1},z_{2})\;\longmapsto\;{1\over1+z_{2}}\pmatrix{w+wz_{2}+
z_{1}^{2}&\sqrt2z_{1}\cr\sqrt2z_{1}&1-z_{2}\cr}\Leqno{OC}$$ defines a
biholomorphic mapping $\phi$ from $U$ to an open subset of $E$ with
$\phi(\6M_{1})\subset\6T$. The inverse $\phi^{-1}$ is given by
$$\pmatrix{x&y\cr y&t\cr}\;\longmapsto\;{1\over1+t}
\Big(x+xt-y^{2},\,\sqrt2y,\,1-t\Big)\,.$$

\bigskip More generally, we consider $\6M:=\6M_{t}$ as a closed
CR-hypersurface of $U_{t}:=\{(w,z_{1},z_{2})\in\CC^{3}:tz_{2}\overline
z_{2}<1\}$ and $\7g=\hol(\6M,0)$ as a Lie algebra of holomorphic vector
fields on $\CC^{3}$. For fixed $t\ne1$
$$N:=\{(w,z_{1},0)\in\CC^{3}:w=z_{1}\overline
z_{1}\}$$ is a
transversal slice in $\6M$ to the infinitesimal action of $\7g$ on
$\6M$. An elementary computation gives
$$\pmatrix{2&2z_{1}\cr2\overline z_{1}&2tz_{1}\overline z_{1}\cr}$$
for the Levi matrix, and hence $4(t-1)z_{1}\overline z_{1}$ for the
Levi determinant at the point $(w,z_{1},0)\in N$. In particular, the
infinitesimal $\7g$-orbit of the origin, that is $\{(w,0,z_{2})\in
\CC^{3}:w\in i\[1\RR,tz_{2}\overline z_{2}<1\}$, coincides with the
set of all Levi degenerate points of $\6M$. Furthermore, the Levi form
at every point of $\6M\backslash N$ is definite if $t>1$ and indefinite
otherwise.

\bigskip The family $\6M_{t}$ yields counter examples to several
statements found in the literature: For instance, contrary to Theorem
4 in \Lit{EBEN}, for every $\6M_{t}$ the vector field
$\xi^{2}=itw^{2}\ddw+itwz_{1}\ddz1 -iz_{1}^{2}\ddz2\in\hol(\6M_{t},0)$
has vanishing $1$-jet at $0$. Furthermore, in \Lit{GAME} it is claimed
that ``$\6M_{1}$ is a new example of a uniformly Levi-degenerate
CR-manifold apart from $\6T$,'' which cannot be true due to Proposition
\ruf{IL}. Finally, translated to our notation, it is claimed in
\Lit{ERSH} \p194 for the model surface $\6M_{0}$ that
``$\dim\hol(\6M,0)\le\dim\hol(\6M_{0},0)=6$ holds for every hypersurface
$\6M$ given by an equation of the form \Ruf{AE}''. But for $\6M_{1}$ this
dimension is $10$.

\KAP{Global}{Global properties}

Let again $V,E,\6C,\6T,\6H$ have the same meaning as in section
\ruf{Preliminaries}. The tube manifold $\6T$ is a homogeneous
$2$-nondegenerate CR-manifold of dimension $5$. In the following we
want to present further CR manifolds of this type that are locally but
not globally CR-equivalent to $\6T$.

Let us denote by $\tilde\6T$ the universal covering of $\6T$ in the
following and by $\mu:\tilde\6T\to\6T$ the corresponding covering map.
For convenience we identify $\tilde\6T$ with $\RR^{2}\times V$ and
$\mu$ with the mapping
$$(r,\phi,v)\;\longmapsto\;e^{r} \Matrix{1+
\cos\phi}{\sin\phi}{\sin\phi}{1-\cos\phi} \;+\;iv\;,\qquad
r,\phi\in\RR,v\in V\;.\Leqno{IS}$$ Clearly, the CR-structure on
$\tilde\6T=\RR^{2}\times V$ is uniquely determined by the condition
that $\mu$ is a local CR-isomorphism.

Our first result states in particular that every continuous
CR-function on $\tilde\6T$ is constant on $\mu$-fibers. For this
denote by $\hat\6T:=\6T\cup \6H$ the convex hull of $\6T$ in $E$.
Also, call for every open subset $U\subset \hat\6T$ a continuous
mapping $f:U\to\CC^{n}$, $n\in\NN$, {\sl holomorphic} if its
restriction to $U\cap \6H$ is holomorphic in the usual sense.

\Proposition{BA} To every continuous CR-function $f$ on
$\tilde\6T$ there exists a unique holomorphic function $h$ on $\hat\6T$
with $f=h\circ \mu$.

\Proof There exists a connected Lie group $G$ (for instance, the group
of all transformations \Ruf{BC} with $\det(g)>0$) acting continuously
on $\hat\6T$ by biholomorphic transformations (in the extended sense
defined above) such that $\6T$ and $\6H$ are $G$-orbits. We may assume
that $G$ is simply connected, otherwise replace $G$ by its universal
covering group. This guarantees that the action of $G$ on $\6T$ lifts
to a transitive action of $G$ on $\tilde\6T$. The sheaf $\5O$ over
$\hat\6T$ of germs of local holomorphic functions (also in the
extended sense) is a Hausdorff space for which the canonical
projection $\pi:\5O\to\hat\6T$ is a local homeomorphism. Now fix an
arbitrary point $a\in\6T$ and choose an open neighbourhood $V$ of $a$
in $\6T$. For $V$ small enough there exists a continuous section
$\sigma:V\to\tilde\6T$, i.e. $\mu\circ\sigma=\id|_{V}$. Since $\6T$ is
the smooth boundary part of the Stein domain $\6H$, the Levi cone at
$a\in\6T$, (which is non zero) points into the direction of $\6H$.
Consequently, there exists a connected open neighbourhood $U$ of $a$
with respect to $\hat\6T$ such that every continuous CR-function on
$V$ has a holomorphic extension to $U$, compare \Lit{BGSS} \p256 or
\Lit{BERO} \p205. Fix now an arbitrary continuous CR-function $f$ on
$\tilde\6T$. For every $g\in G$ there exists a unique holomorphic
function $f_{g}$ on $U$ which coincides with $f\circ g\circ\sigma$ on
$U\cap V$. Denote by $\5F_{g}\subset\5O$ the subset of all germs
induced by the holomorphic function $f_{g}\circ g^{-1}$, defined on
$g(U)$. Clearly, every $\5F_{g}$ is connected. Observe that for
sufficiently close $g_{1},g_{2}\in G$, the intersection
$\5F_{g_{1}}\cap\5F_{g_{2}}$ is not empty. Since $G$ is connected,
there is a unique connected component $\5F$ of $\5O$ containing all
$\5F_{g}$, $g\in G$. We claim that $\pi:\5F\to\hat\6T$ is a covering
map. Since $\pi$ is a local homeomorphism it is sufficient to show
that every continuous curve $\gamma:[0,1]\to\hat\6T$ has a lifting to
a continuous curve $\tau:[0,1]\to\5F$ with
$\pi\circ\tau=\gamma$. Without lost of generality we may consider only
curves $\gamma$ in $\6H$. The claim follows from the existence of a
continuous curve $t\mapsto g_{t}$ in $G$ with
$\gamma(t)=g_{t}(\gamma(0))$ for all $t\in[0,1]$. Since $\hat\6T$ is
simply connected we get that $\pi:\5F\to\hat\6T$ is a
bijection. Therefore $\5F$ is a continuous section in $\5O$ over
$\hat\6T$ and determines the required holomorphic function $h$ on
$\hat\6T$.\qed

\Corollary{} To every $\tilde g\in\Aut(\tilde\6T)$ there exists a
unique $g\in\Aut(\6T)$ with ${\mu\circ\tilde g=g\circ\mu}$. The mapping
$\Aut(\tilde\6T)\to\Aut(\6T)$, $\,\tilde g\mapsto g$, realizes
$\Aut(\tilde\6T)$ as universal covering group of $\Aut(\6T)$. In
particular, $\Aut(\tilde\6T)$ is a simply connected Lie group of
dimension $7$ with two connected components, acting transitively on
$\tilde\6T$.

\Proof Denote by $\iota:\6T\hookrightarrow E$ the canonical
injection. By Proposition \ruf{BA} the mapping
${\iota\circ\mu\circ\tilde g}$ is constant on $\mu$-fibers and hence
factors over $\mu$. The group $\Aut(\6T)$ is given by all affine
transformations \Ruf{BC} and hence there is a canonical isomorphism of
Lie groups
$$\Aut(\6T)\;\cong\;\GL(\6C)\ltimes V\,,$$ where
$$\GL(\6C):=\{g\in\GL(V):g(\6C)=\6C\} \;\cong\;(\GL(2,\RR)/\{\pm e\})\,,$$
$e\in\GL(2,\RR)$ is the unit matrix and the semi-direct product refers
to the canonical injection $\rho:\GL(\6C)\hookrightarrow\GL(V)$. In
particular, $\Aut(\tilde\6T)\cong\tilde{\GL}(\6C)\ltimes V$ with
$\tilde{\GL}(\6C)\cong\tilde{\GL}(2,\RR)$ the universal covering group
of $\GL(\6C)$.\qed

\Ues{A CR-deformation family for $\9{\6T}\!$} ~For every $\psi\in\RR$
and $\mu_{\psi}:=
\Matrix{\cos\psi/2}{\sin\psi/2}{-\sin\psi/2}{\cos\psi/2}$ define
$\lambda_{\psi}\in\GL(\6C)$ by
$\lambda_{\psi}v=\mu_{\psi}v\mu_{\psi}'$. For every
$z=(s,\psi,w)\in\RR^{2}\Times V$ then
$$\theta_{z}(r,\phi,v):=(r+s,\[1\phi+\psi,
\[1e^{s}\lambda_{\psi}v+w)\qquad \Leqno{UO}$$ defines an affine
transformation $\theta_{z}\in\Aut(\tilde\6T)$ and
$\Theta:=\{\theta_{z}:z\in\RR^{2}\Times V\}$ is a Lie group acting
freely and transitively on $\tilde\6T$. Clearly, $\Theta$ has many
discrete subgroups $\Gamma$, each of which gives a CR-manifold
$\tilde\6T/\Gamma$ locally CR-equivalent to $\6T$. Here we restrict
our attention to the following one-parameter family
$(\Gamma_{t})_{t>0}$ of discrete subgroups: For every real $t>0$ put
$\gamma_{t}:=\theta_{(0,2\pi t,0)}$,
$\;\Gamma_{t}:=\{\gamma_{nt}:n\in\ZZ\}$ and
$\6T_{t}:=\tilde\6T/\Gamma_{t}$. Then every $\6T_{t}$ is isomorphic to
$\6T$ as real-analytic manifold while $\6T_{1}$ is equivalent to $\6T$
as CR-manifold. We will see later that actually the $\6T_{t}$ are
pairwise nonequivalent as CR-manifolds, compare Proposition
\ruf{PT}. The family $(\6T_{t})_{t>0}$ may be considered as a
CR-deformation family of $\6T\cong\6T_{1}$: For
$\RR^{+}:=\{t\in\RR:t>0\}$ consider the CR-manifold
$\tilde\6T\Times\RR^{+}$, on which $\ZZ$ acts freely by
$(z,t)\mapsto(\gamma_{nt}z,t)$, $n\in\ZZ$. Then
$(\tilde\6T\Times\RR^{+})/\ZZ$ is a real-analytic CR-manifold, and the
canonical projection $\pi:(\tilde\6T\Times\RR^{+})/\ZZ\to\RR^{+}$ is a
CR-mapping whose fibers give the family $(\6T_{t})_{t>0}$.

\Proposition{IU} For every $t>0$ and every $a\in\6T_{t}$ the
cardinality $\sigma_{t}(a)$ of the set
$$\Sigma_{t}(a):=\{z\in\6T_{t}:f(z)=f(a)\steil{for every continuous
CR-function $f$ on}\6T_{t}\}\Leqno{YW}$$ is given by
$$\sigma_{t}(a)=\cases{p&$t=p/q$\steil{for relatively prime
integers}$p,q>0$\cr\noalign{\vskip2pt}
\infty&{\rm in all other cases.}\cr
}$$

\Proof The group $\Gamma:=\Gamma_{1}/(\Gamma_{1}\cap\Gamma_{t})$ acts
freely on $\6T_{t}$, and by Proposition \ruf{BA} every continuous
CR-function on $\6T_{t}$ is constant on $\Gamma$-orbits. In case $t$
irrational $\Gamma\cong\ZZ$ implies $\sigma_{t}(a)=\infty$. Therefore
we may assume in the following that $t=p/q$ with relatively prime
integers $p,q>0$. Since then $\Gamma$ has order $p$ we get
$\sigma_{t}(a)\ge p$. On the other hand, the CR-manifold $\6T_{1/q}$
is separable with respect to real-analytic CR-functions. Indeed, the
orbits of the finite subgroup $\Lambda:={\{\lambda^{}_{2\pi s}:s\in
q^{-1}\ZZ\}}\subset\GL(V)\subset\GL(E)$ are separated by the
$\Lambda$-invariant holomorphic polynomials on $E$. Since we have a
CR-covering map $\6T_{p/q}\to\6T_{1/q}$ of degree $p$ this implies
${\sigma_{t}(a)\le p}$.\phantom{QQQ}\qed

It is not difficult to see that the centralizer of $\gamma_{t}$ in
$\Aut(\tilde\6T)$ acts transitively on $\tilde\6T$ if and only if
$t\in\NN$. As a consequence, $\6T_{t}$ is homogeneous as CR-manifold
if and only if $t$ is an integer. In any case, the centralizer of
$\gamma_{t}$ contains the 1-parameter subgroup
$\{\gamma_{s}:s\in\RR\}$ of $\Aut(\tilde\6T)$, which therefore also
acts on $\6T$.

\Proposition{PZ} For every irrational $t>0$ the following properties
hold: \0 For every $a\in\6T_{t}$, the set $\Sigma_{t}(a)$ defined in
\Ruf{YW} is the circle $\Sigma_{t}(a)=\{\gamma_{s}(a):0\le
s<t\}$. \1 There exists a closed hypersurface
$\6N_{t}\subset\6T_{t}$ such that every continuous CR-function on
$\6T_{t}$ is real-analytic on the complement
$\6T_{t}\backslash\6N_{t}$.\1 There exists a continuous CR-function on
$\6T_{t}$ which is not real-analytic.

\Proof Fix a continuous CR-function $f$ on $\6T_{t}$ and denote by
$\mu:\tilde\6T\to\6T$, $\mu_{t}:\tilde\6T\to\6T_{t}$ the canonical
projections. By Proposition \ruf{BA} there exists a holomorphic
function $h$ on $\6H\cup\6T$ with $f\circ\mu_{t}=h\circ\mu$. Since,
for every $a\in\6T_{t}$, the function $f$ is constant on a dense
subset of the circle $\{\gamma_{s}(a):s\in\RR\}$ we get
$\Sigma_{t}(a)\subset\{\gamma_{s}(a):s\in\RR\}$. This implies that
$f\circ\mu_{t}$ is $\gamma_{s}$-invariant for every $s\in\RR$ and
hence that $h$ is invariant under $\SO(2)$ acting on $E$ as in
\Ruf{BC}. Now the proposition follows from the next lemmata \ruf{QT}
and \ruf{QR}, which may be of independent interest. The inclusion
$\{\gamma_{s}(a):s\in\RR\}\subset\Sigma_{t}(a)$ follows from the fact
that the mapping $\phi$ defined in Lemma \ruf{QT} separates the
$\SO(2)$-orbits in $\6T$.\qed

\Lemma{QT} Consider in $\6T$ the hypersurface
$$\6N:=\{c+i(rc+se):c\in\6C\steil{and}r,s\in\RR\},\quad e\hbox{~\rm
the unit matrix,}$$ and denote by $\6W:=\6T\backslash\6N$ its
complement. Let $\tr$ be the normalized trace on $E$
(i.e. $\tr(e)=1$) and define the $\SO(2)$-invariant holomorphic
mapping $\phi:E\to\CC^{2}$ by
$\phi(z):=\big(\tr(z),\tr(z^{2})-\tr(z)^{2}\big)$. Then
$$\eqalign{\phi(\6W)\;&=\;\big\{(x_{1}+iy_{1},x_{2}+iy_{2})\in\CC^{2}:
x_{1}>0,\; x_{2}<x_{1}^{2}-(2x_{1})^{-2}y_{2}^{2}\big\}\cr
\phi(\6N)\;&=\;\big\{(x_{1}+iy_{1},x_{2}+iy_{2})\in\CC^{2}: x_{1}>0,\;
x_{2}=x_{1}^{2}-(2x_{1})^{-2}y_{2}^{2}\big\}\,,\cr }$$ and to every
continuous $\SO(2)$-invariant CR-function $h$ on $\6T$ there exists a
unique continuous function $f$ on $\phi(\6T)$ with $h=f\circ\phi$ such
that the restriction of $f$ to the interior $\phi(\6W)$ is
holomorphic. In particular, every continuous $\SO(2)$-invariant
CR-function on $\6T$ is real-analytic on the dense domain
$\6W\subset\6T$.

\Proof The expressions for $\phi(\6W)$ and $\phi(\6N)$ are checked by
direct computation. Every continuous $\SO(2)$-invariant function $h$ on
$\6T$ is of the form $h=f\circ\phi$ with $f$ a continuous function on
$\phi(\6T)$. In case $h$ is CR in addition, $h$ admits an extension to
an $\SO(2)$-invariant holomorphic function on $\6H$. This implies that
$f$ is holomorphic on $\phi(\6W)=\phi(\6H)$.\qed

For the proof of \ruf{PZ}.(iii) it would be enough to find directly a
continuous function $g$ on $\phi(\6T)$ such that the pull-back
$f:=g\circ\phi$ is CR but not real-analytic on $\6T$. Instead, we use
invariant integration for the explicit determination of such an $f$.

\Lemma{QR} There exists a continuous $\SO(2)$-invariant function on
$\6T$ that is not real-analytic.

\Proof Let $\sqrt s$ be the unique continuous branch of the square
root on $H:=\{s\in\CC:\Re(s)\ge0)\}$ with $\sqrt1=1$. Then
$\lambda(z)=z_{11}$ defines a complex linear form $\lambda:E\to\CC$
with $\lambda(\hat\6T)\subset H$ for $\hat\6T:=\6T\cup\6H$. The
function $\phi(z):=\sqrt{\lambda(z)}$ is continuous on $\hat\6T$ and
holomorphic on $\6H$. For easier computation let us introduce on $E$
the coordinates $(w,z_{1},z_{2})\mapsto
\Matrix{w-z_{1}}{z_{2}}{z_{2}}{w+z_{1}}$, and then
$\lambda(w,z_{1},z_{2})=w-z_{1}$. For every $t\in\RR$ define
$\gamma_{t}\in\GL(E)$ by $$(w,z_{1},z_{2})\;\longmapsto\;(w,z_{1}\cos
t+z_{2}\sin t,-z_{1}\sin t+z_{2}\cos t)\,.$$ Then
$$h(z):=\intop^{2\pi}_{0}\!\phi(\gamma_{t}z)\,dt$$ defines a
continuous $\SO(2)$-invariant function $h$ on $\hat\6T$ that is
holomorphic on $\6H$. As local uniform limit of holomorphic functions
therefore the restriction $f:=h|_{\6T}$ is a continuous
$\SO(2)$-invariant CR-function on $\6T$. Define on the open unit disk
$\Delta\subset\CC$ the holomorphic function $g$ by
$$g(s):=h(1,s,0)=\intop_{0}^{2\pi}\sqrt{1-s\cos t}\,dt=
\sum_{k=0}^{\infty}{1/2\choose k}c^{}_{k}s^{k},\quad
c^{}_k:=\intop_{0}^{2\pi}(\cos t)^{k}dt\,,$$ and denote by $R\ge1$ the
radius of convergence of the power series expansion. From the
recursion $c_{0}=2\pi$, $c_{1}=0$ and $c^{}_{k}={k-1\over
k^{}}c^{}_{k-2}$ for $k\ge2$ we get $R=1$. Because of $(1,s,0)\in\6T$
for $|s|=1$ the function $g$ has a continuous extension to the closure
$\overline\Delta$.

\noindent Now assume that $f$ is real-analytic on $\6T$. Then $g$
would have a holomorphic extension to an open neighbourhood of
$\overline\Delta$ in $\CC$. This contradicts $R=1$, that is, $f$
cannot be real-analytic on $\6T$.\qed

\medskip

\Proposition{PT} The manifolds $\6T_{s},\6T_{t}$ are globally
CR-equivalent if and only if $s=t$.

\Proof Assume that $\6T_{s},\6T_{t}$ are globally CR-equivalent. Then
there exists a transformation $g\in\Aut(\tilde\6T)$ with
$g\[2\Gamma_{\!s}=\Gamma_{\!t}\[1g$, that is,
$g\[2\gamma_{s}=\gamma_{r}\[1g$ with $r=\pm t$. In case $s\in\NN$ the
CR-manifold $\6T_{s}$ is homogeneous. But then
$\sigma_{t}(a)=\sigma_{s}(b)$ for all $a\in\6T_{a}$ and $b\in\6T_{s}$
implies $s=t$ as a consequence of Proposition \ruf{IU}. We
therefore may assume $s,t\notin\NN$. Inspecting the
homomorphic image of the equation $g\[2\gamma_{s}=\gamma_{r}\[1g$ in
$\Aut(\6T)$ shows $g\in\Theta$ and thus $s=t$, see \Ruf{UO} for the
definition of the group $\Theta$.\qed

\KAP{bounded}{The bounded realization}

Let again $E=\{z\in\CC^{2\times2}:z'=z\}$ be the linear space of
all complex symmetric $2\times2$-matrices with $e=\One_{2}\in E$ the
unit matrix. Put $$\eqalign{\6D:&=\{z\in E:e{-}z\overline z\hbox{
positive definite}\}\cr \6R:&=\{z\in E:\det(e{-}z\overline
z)=0,\;\tr(e{-}z\overline z)>0\} \;=\;\{z\in\partial \6D:z\overline
z\ne e\}\cr\6S:&=\{z\in E:z\hbox{ unitary}\}\;\cong\;
\U(2)/\O(2)\;\cong\;(S^{1}\times S^{2})/\ZZ_{2}\;,\cr}\Leqno{UV}$$
where the action of $\ZZ_{2}$ is generated by the antipodal map in
each of the two spheres $S^{1}, S^{2}$. Then $\6D$ is the open unit
ball in $E$ with respect to the operator norm and is also called a
{\sl Lie ball}. The boundary $\partial \6D$ decomposes into its smooth
part $\6R$ and the Shilov boundary $\6S$, which is a totally real
submanifold of $E$.

The CR-manifold $\6R$ is fibered in its {\sl holomorphic arc
components,} which all are affinely equivalent to the open unit disk
$\Delta\subset\CC$, compare \Lit{WOLF}. By definition, $A\subset\6R$
is a holomorphic arc component of $\6R$, if it is minimal with respect
to the property: {\sl $A\ne\emptyset$, and $f(\Delta)\subset A$ for
every holomorphic map $f:\Delta\to E$ with $f(\Delta)\subset\6R$ and
$f(\Delta)\cap A\ne\emptyset$.} For instance, the set $A$ of all
diagonal matrices $z\in E$ with $z_{11}=1>|z_{22}|$ is such an arc
component, and every other is of the form $uAu'$ with
$u\in\U(2)$. This implies that the space of all arc components in
$\6R$ can be identified with
$\U(2)/\big(\!\O(1)\times\U(1)\big)\cong\PP_{3}(\RR)$. Since every
$g\in\Aut(\6R)$ respects holomorphic arc components there is an
$\Aut(\6R)$-equivariant fiber bundle $\Xi:\6R\to\PP_{3}(\RR)$ with
fibers the holomorphic arc components of $\6R$.

It is well known that the {\sl Cayley transformation} $z\mapsto
(z-e)(z+e)^{-1}$ defines a biholomorphic map $\gamma:\6H\to \6D$ and
also gives a CR-isomorphism $\6T\buildrel\simeq\over\to\{z\in
\6R:\det(e-z)\ne0\}$. In fact, $E$ can be considered as Zariski-open
subset of a nonsingular quadric $Z\subset\PP_{4}(\CC)$ in such a way
that $\gamma$ extends to a biholomorphic automorphism of $Z$, and then
$\gamma^{2}(z)=-z^{-1}$ for all invertible $z\in E$. Also, there
exists a unique antiholomorphic involution $\tau:Z\to Z$ with
$\Fix(\tau)=\6S$, and then
$$\eqalign{G:&=\{g\in\Aut(Z):g(\6D)=\6D\}=\{g\in\Aut(Z):g(\6R)=\6R\}\cr
&=\{g\in\Aut(Z):\tau g\tau=g\}^{0} \;\cong\;\SO(2,3)^{0}\cr }$$ acts
transitively on $\6D$, $\6R$, $\6S$. Via these actions the groups $G$,
$\,\Aut(\6D)$ and $\Aut(\6R)$ are canonically isomorphic and hence are
identified in the following. We also identify the Lie algebra
$$\7g:=\aut(\6D)=\big\{(a+cz+zc'-z\overline az)\ddz{}:a\in
E,\,c\in\7u(2)\big\}$$ of $G$ with $\aut(\6R)=\hol(\6R)$, where
$\7u(2)\subset\7{gl}(2,\CC)$ is the Lie subalgebra of all skew
hermitian matrices.

The complexification $\7l=\7g\oplus i\7g$ has the decomposition
$$\7l=\bigoplus_{\nu\in\ZZ^{2}}\7l^{v}\Steil{with}\7l^{\nu}:=
\{\xi\in\7l:[\zeta_{j},\xi]=\nu_{j}\[1\xi\steil{for}j=1,2\}$$ with
$\zeta_{j}\in\7l$ defined by \Ruf{FX} and \Ruf{DR} again being a
choice of root vectors $\xi^{\nu}\in\7l^{\nu}$. But unlike to \Ruf{OR}
the real form $\7g$ is embedded in a different way here:
$$\7g\cap(\7l^{\nu}+\7l^{-\nu})=\cases{\RR(\xi^{\nu}+\xi^{-\nu})
\oplus\RR(i\xi^{\nu}-i\xi^{-\nu})&$\nu\ne0$ \cr
\RR i\zeta_{1}\oplus\RR i\zeta_{2}&otherwise$\;.$\cr}\Leqno{TR}$$

\medskip $\6R$ has the remarkable property that every CR-equivalence
between domains $\6U,\6V$ in $\6R$ extends to a CR-automorphism of
$\6R$, compare \Lit{KAZT}. In particular, $\Aut(\6U)\cong\{g\in
G:g(\6U)=\6U\}$ and the groups $\Aut(\6U)$, $\,\Aut(\6V)$ are
conjugate in $\Aut(\6R)$ for any pair of CR-equivalent domains
$\6U,\6V\subset\6R$. Via the Cayley transformation $\6R$ contains a
copy of the CR-manifold $\6T$ and hence may be considered as a
CR-extension of $\6T$. We may ask whether $\6R$ is maximal with
respect to CR-extensions, that is, whether $\6R$ can be a proper
domain in some other connected CR-manifold. A partial answer is given
by the following result.

\Proposition{} Let $M$ be a connected, locally homogeneous CR-manifold
and let $D\subset M$ be a domain that is CR-isomorphic to a covering of
$\6R$. Then $D=M$.

\Proof Assume to the contrary $D\ne M$ and fix a boundary point $c$ of
$D$ in $M$. Since $M$ is locally homogeneous there exists a connected
open neighbourhood $U\subset M$ of $c$ and a vector field
$\eta\in\hol(U)$ having an integral curve $\gamma:[0,1]\to U$ with
$\gamma(0)=c$ and $a:=\gamma(1)\in D$. Since $D$ is CR-equivalent to a
covering of $\6R$ there exists a vector field $\xi\in\hol(D)$ having
the same germ as $\eta$ at $a$. Since $\xi$ is complete on $D$ there
exists an integral curve $g:\RR\to D$ of $\xi$ with $g(0)=a$. But
then, for every $s\in[0,1]$, we have $g_{-s}=\gamma_{1-s}\in D$. In
particular, $c=g_{-1}\in D$ gives a contradiction.\qed

\Proposition{UG} Let $M$ be a homogeneous, connected, simply connected
real-analytic CR-manifold which is locally CR-isomorphic to $\6R$. Then
$M$ is CR-equivalent to the universal covering $\tilde\6U$ of a
homogeneous domain $\6U\subset\6R$.

\Proof For every $a\in M$ let $\5F_{a}$ be the set of all germs at
$a$ of CR-equivalences $U\to \6V$, where $U$ is a an open
neighbourhood of $a\in M$ and $\6V$ is open in $\6R$. Consider the
disjoint union $\5F$ of all $\5F_{a}$, $a\in M$, in the usual way as
sheaf over $M$ with sheaf projection $\mu:\5F\to M$ and fix a
connected component $N$ of $\5F$. By assumption there is a connected
Lie group $H$ of CR-automorphisms of $M$ acting transitively on
$M$. This implies that $\mu:N\to\6M$ is a covering map and hence
bijective. The evaluation map $N\to\6R$ defines a CR-covering map
onto the image $\6U\subset\6R$. Since this map is $H$-equivariant,
$\6U$ is a homogeneous domain in $\6R$.\qed

Proposition \ruf{UG} reduces the classification problem for
CR-manifolds locally CR-isomorphic to $\6R$ to the study of
homogeneous domains in $\6R$. But these are in 1-1-correspondence to
$G$-homogeneous domains in $\PP_{3}(\RR)$. More precisely, every
homogeneous domain $\6U\subset\6R$ is the full $\Xi$-pre-image of
$U:=\Xi(\6U)$. Indeed, fix an arbitrary holomorphic arc component $A$
of $\6R$ with $A\cap\6U\ne\emptyset$. Then the group
$H:=\{g\in\Aut({\6U}):g(A)=A\}$ has the orbit $A\cap\6U$ that is open
in $A\cong\Delta$. Therefore $H$ acts transitively on $A$ implying the
claim. As a consequence, $\6U\to U$ is a fiber bundle with
contractible fiber, implying that the homogeneous domain
$\6U\subset\6R$ and its $\Xi$-image $U\subset\PP_{3}(\RR)$ always have
the same homotopy type. In particular, $\6R$ has fundamental group
$\ZZ_{2}$.

For certain homogeneous domains in $\6R$ the same argument as in the
proof of Proposition \ruf{BA} yields the following holomorphic
extension property.

\Proposition{TZ} Let $\mu:\tilde{\6U}\to\6U$ be the universal covering
of a homogeneous domain $\6U\subset\6R$ and assume that
$\Aut(\6U)\subset\Aut(\6R)$ acts transitively on $\6D$. Then to every
continuous CR-function $f$ on $\tilde{\6U}$ there exists a unique
continuous function $h$ on $\6U\cup\6D$ with the following two
properties. \0$f=h\circ \mu$, \1 $h$ is holomorphic on $\6D$.\Formend

\bigskip\noindent 
To get further examples of homogeneous domains in
$\6R$ write $$G_{Y}:=\{g\in
G:g(Y)=Y\}\steil{and}\7g_{Y}:=\{\xi\in\7g:\exp(\RR\xi)\subset
G_{Y}\}$$ for every subset $Y\subset Z$. Then every subgroup $G_{a}$,
$a\in\6S$, has a unique open orbit in $\6R$. In case $a=-e$ this orbit
is just the image $\gamma(\6T)\subset\6R$ of the tube manifold $\6T$.
Also, for every holomorphic arc component $A$ of $\6R$ the group
$G_{A}$ has a unique open orbit in $\6R$. Even the intersection
$G_{A}\cap G_{a}$ has a unique open orbit in $\6R$, provided $a\in\6S$
is in the closure of $A$ with respect to $\partial\6D$. All these
groups act transitively on $\6D$ and thus give homogeneous domains in
$\6R$ with the holomorphic extension property of Proposition \ruf{TZ}.
But there also exist homogeneous domains $\6U\subset\6R$ such that
$\Aut(\6U)$ is not transitive on $\6D$:

\Example{} Let $F\subset E$ be the $\CC$-linear subspace
of all diagonal matrices. Then the intersection $F\cap \6D$ is
biholomorphically equivalent to the bidisk $\Delta^{2}\subset\CC^{2}$
and $G_{F\cap \6D}$ has Lie algebra
$$\7g_{F\cap \6D}:=\big\{(a+cz+zc'-z\overline az)\ddz{}:a\in
F,\,c\in F\cap\7u(2)\big\}\;\cong\;
\7{sl}(2,\RR)\times\7{sl}(2,\RR)\,.$$ The corresponding open orbit
$\6U\subset\6R$ has the following equivalent description, which may be
of interest for itself: The group $H:=\SL(2,\RR)\times\SL(2,\RR)$ acts
on the complex manifold $X:=\SL(2,\CC)$ by $z\mapsto gzh^{-1}$. For
every $\Matrix abcd\in X$ put $\delta(z):=\det(z{+}\overline
z)-2=2\Re\big(a\overline d-b\overline c\big)$. Then $\delta$ is an
$H$-invariant function on $X$ with $\delta(X)=\RR$. It is easily
checked that multiplication by $j:=\Matrix0ii0$ induces a
biholomorphic automorphism of $X$ with $\delta(jz)=-\delta(z)$. The
critical set of $\delta$ is $\{z\in X:\overline z=\pm
z\}=\SL(2,\RR)\cup j\cd\SL(2,\RR)$. In particular, $\delta$ has
critical values $\pm2$. For every $\epsilon\in\{\pm i,0\}$ consider
the orbit $M_{\epsilon}:=H\Matrix1\epsilon01$. The inversion $z\mapsto
z^{-1}$ maps $M_{\epsilon}$ to $M_{-\epsilon}$ and $\overline
{M_{i}}\cap\overline{M_{-i}}=M_{0}=\SL(2,\RR)$. It can be shown that
$M_{i}$ is CR-equivalent to $\6U\subset\6R$. Furthermore, $M_{i}$ is
diffeomorphic to $\6C\times\SL(2,\RR)$ and hence has fundamental group
$\ZZ^{2}$ (recall that the tube manifold $\6T$ is diffeomorphic to
$\6C\times V$).

\Example{} The Lie group $\SU(2)/\{\pm e\}\cong\SO(3)$ acts on $E$ by
$z\mapsto gzg'$. Choose an $\RR$-linear subspace $W$ of $E$ that is
invariant under this action, e.g. $W=\{z\in E:z_{11}-\overline
z_{22}=z_{12}+\overline z_{12}=0\}$. Denote by $\tau:E\to E$ the
unique conjugate linear involution fixing every point of $W$. Then
$\tau(\6D)=\6D$ and the intersection $W\cap \6D$ is a euclidian ball
in $W$. The group $G_{W\cap \6D}=\{g\in G:\tau g\tau=g\}$ has Lie
algebra
$$\7g_{W\cap \6D}=\big\{(a+cz+zc'-z\overline az)\ddz{}:a\in
W,\,c\in\7{su}(2)\big\}\;\cong\;\7{sl}(2,\CC)\,.$$ Again, the
corresponding open orbit $\6U\subset\6R$ has an equivalent
realization: Let $X:=\SL(2,\CC)$ and consider $H:=\SL(2,\CC)$ as real
Lie group. Then $H$ acts on $X$ by $z\mapsto gzg^{*}$, where $g^{*}$
is the conjugate transpose of $g$. Clearly, $z\mapsto -z$ commutes
with the action of $H$ and hence permutes $H$-orbits. The function
$\delta(z):=\det(z{+}z^{*})-2=a\overline d+d\overline a-b\overline
b-c\overline c$ on $X$ is $H$-invariant with critical values $\pm2$
and critical set $\{z\in X:z^{*}=\pm z\}$. This set is the union of
three totally real orbits, more precisely:\nline {\sl Critical value
$2$:} There are the two critical orbits $N:=\{h\in X:h^{*}=h>0\}$ and
$-N=\{h\in X:h^{*}=h<0\}$. A diffeomorphism
$\{h=h^{*}\in\CC^{2\times2}:\tr(h)=0\}\to N$ is given by the
exponential mapping. By elementary calculus it is seen that $\delta$
attains a local maximum at every point of $\pm N$. As a consequence,
$N$ is the only $H$-orbit in $X$ having $N$ in its closure (and the
same with $-N$).\nline {\sl Critical value $-2$:} The only critical
orbit is $N:=\{ih:h=h^{*}\in\CC^{2\times2},\,\det(h)=-1\}$. The
hypersurface orbit $M:=H \Matrix1ii0 $ has $N$ in its closure. It can
be shown that $M$ is CR-equivalent to $\6U$.

\KAP{final}{Some final remarks.} 

The idea of proof in Propositions \ruf{BA} and \ruf{TZ} can also be
applied in higher dimensions: For $n\ge3$ fixed let
$V\subset\RR^{n\Times n}$ be the subspace of all symmetric $n\Times
n$-matrices and $E:=V\oplus iV\subset\CC^{n\Times n}$ its
complexification. For all $p,q\ge0$ with $k:=n-(p+q)\ge0$ denote by
$\6C_{p,q}\subset V$ the cone of all real symmetric $n\Times
n$-matrices of type $(p,q)$, that is, having $p$ positive and $q$
negative eigenvalues. The group $\GL(n,\RR)$ acts on $V$ by $x\mapsto
gxg'$, and the cones $\6C_{p,q}$ are the corresponding
orbits. Because of $\6C_{q,p}=-\6C_{p,q}$ we may restrict our
attention to the special case $p\ge q$. The cone $\Omega:=\6C_{n,0}$
is convex open, and $\6H:=\Omega\oplus iV\subset E$ is a symmetric
tube domain (Siegel's upper half plane of rank $n$ up to the
factor~$i$). The group $\Aut(\6H)$ consists of all biholomorphic
transformations \Ruf{WC} with $\Matrix abcd$ in the symplectic
subgroup $\Sp(n,\RR)\subset\SL(2n,\RR)$, compare \Lit{KLIN} \p345.

In analogy to our setting in Section \ruf{tube} with $n=2$ we denote
for every $p,q$ by $\6T_{p,q}:=\6C_{p,q}\oplus iV$ the tube manifold
over the cone $\6C_{p,q}$. Then $\6T_{p,q}$ is open precisely if $k=0$
and is closed if $k=n$. In all other cases $\6T_{p,q}$ is a
$2$-nondegenerate CR-manifold and $\Aut(\6T_{p,q})$ consists of all
transformations \Ruf{BC} (with $\GL(2,\RR)$ replaced by $\GL(n,\RR)$,
compare \Lit{KAPP}). In \Lit{KAPP} also the following holomorphic
extension property has been shown: {\sl In case $pq\ne0$ every
continuous CR-function on $\6T_{p,q}$ has a holomorphic extension to
all of $E$. In case $p>0$ every continuous CR-function on $\6T_{p,0}$
has a holomorphic extension to $\6H$ that is continuous up to
$\6T_{p,0}\subset\overline\6H$ in a certain sense.}

Note that every $\6T_{p,q}$ is homogeneous under $\GL(n,\RR)^0\ltimes
V.$ Applying the decomposition theorem of Mostow one deduces that
$\6T_{p,q}$ is a bundle with contractible fibers over the homogeneous
space $\SO(n)/S\big(\O(p)\Times\O(q)\Times\O(k)\big)$. Explicit
computation shows that $\SO(n)$ is a maximal compact subgroup of
$\GL(n,\RR)^{0}$ and that for the diagonal matrix
$a:=\One_{p}\Times{-}\One_{q}\Times0_{k}\in\6T_{p,q}$ the subgroup
$S\big(\O(p)\Times\O(q)\Times\O(k)\big)\subset\SO(n)$ is maximal
compact in the isotropy subgroup $\{g\in\GL(n,\RR^{n})^{0}:gag'=a\}$.
As a consequence, employing the long exact homotopy sequence, we
determine the fundamental group of $\6T_{p,q}$ for $n\ge3$ and
$0<p+q<n$ as
$$\pi_{1}(\6T_{p,q})=\cases{ Q_{8}&$n=3,\,p=q=1$\cr
\ZZ_{2}\oplus\ZZ_{2}&$n>3,\,p>0,\,q>0$\cr\ZZ_{2}&otherwise$\;,$\cr}$$
where $Q_{8}$ is the quaternion group of order $8$.

For the universal covering $\mu:\tilde\6T_{p,q}\to\6T_{p,q}$ of
$\6T_{p,q}$ we get the following higher dimensional analog of
Proposition \ruf{BA}.

\Proposition{OD} In case $0<p+q<n$ every continuous CR-function on
$\tilde\6T_{p,q}$ is constant on $\mu$-fibers.

\Proof Let $H$ be the group of all transformations $z\mapsto
gzg'+iv$ on $E$ with $g\in\GL(n,\RR)$ and $v\in V$. To begin with
assume first that $p\cd q>0$ and fix a continuous CR-function $f$ on
$\tilde\6T_{p,q}$. Let $U\subset E$ be the smallest open $H$-invariant
subset containing $ \6T_{p,q}$, that is, $U$ is the union of all
$\6T_{p',q'}$ with $p'\ge p$ and $q'\ge q$. Denote by $\5O$ the sheaf
over $U$ of all germs of holomorphic functions. Since for every
$a\in\6T_{p,q}$ the Levi cone at $a$ spans the full normal space to
$\6T_{p,q}$ in $E$, for every $x\in\tilde\6T_{p,q}$ the function $f$
induces a germ $f_{x}\in\5O_{\mu(x)}$. Denote by $\5F$ the connected
component of $\5O$ containing all $f_{x}$ with $x\in\tilde\6T_{p,q}$.
As in the proof of Proposition \ruf{BA} it follows that $\5F$ is a
covering of $U$. As a consequence of the following Lemma \ruf{LE}
every loop in $\6T_{p,q}$ is zero-homotopic in $U$. Therefore
$\{f_{x}:x\in\tilde\6T_{p,q}\}$ is a trivial covering of $\6T_{p,q}$,
that is, $f=g\circ\mu$ for some CR-function on $\6T_{p,q}$.\nline It
remains to consider the case $q=0$. But then $\6T_{p,0}\cup\6H$ is
simply connected and the claim follows in the same way as in the proof
of Proposition \ruf{BA} (for every open subset $V$ of
$\6T_{p,0}\cup\6H$ a function $g$ on $V$ should be called {\sl
holomorphic} if its restriction to $V\cap\6H$ is holomorphic in the
usual sense and if for every $a\in V\cap\6T_{p,q}$ there is a
neighbourhood $W$ of $a$ in $\6T_{p,q}$ and a wedge $\Gamma$ with edge
$W$ such that the restriction $g_{|\Gamma}$ extends continuously to
$W$).

\Lemma{LE} For every $p,q$ with $p+q<n$ the canonical injection
$\6C_{p,q}\hookrightarrow\6C_{p+1,q}\cup\6C_{p,q}\cup\6C_{p,q+1}$
induces the trivial homomorphism between the corresponding fundamental
groups.

\Proof Consider the 1-parameter subgroup of $\SO(n)$ given by
$$g_{t}:=\pmatrix{\cos t&-\sin t&\cr\sin t&\cos t&\cr&&\One_{n{-2}}}$$
for all $t\in\RR$. Furthermore fix diagonal matrices $a,b\in\6C_{p,q}$
with $a_{11}=1=-b_{11}$ and $a_{22}=b_{22}=0$. Then the loops
$\gamma,\sigma:[0,\pi]\to\6C_{p,q}$ defined by
$\gamma(t):=g^{}_{t}\]2ag'_{t}$ and
$\sigma(t):=g^{}_{t}\]2bg'_{t}$ generate the group
$\pi_{1}(\6C_{p,q})$. For every $s\in\RR$ let $a(s),b(s)$ be the
matrices obtained from $a,b$ by putting $a(s)_{22}=b(s)_{22}=s$ and
leaving all other entries unchanged. Then $a(0)=a$ and
$a(s)\in\6C_{p+1,q}$ for all $s>0$. But then the loops
$\gamma_{s}:[0,\pi]\to\6C_{p,q}\cup\6C_{p+1,q}$, $s\ge0$, defined by
$\gamma_{s}(t):=g^{}_{t}\]2a(s)g'_{t}$ give a homotopy from
$\gamma$ to the constant loop $\gamma_{1}$. In the same way the loops
$\sigma_{s}:[0,\pi]\to\6C_{p,q}\cup\6C_{p,q+1}$, $s\le0$, defined by
$\sigma_{s}(t):=g^{}_{t}\]2b(s)g'_{t}$ give a homotopy from
$\sigma$ to the constant loop $\sigma_{-1}$.\qed

\medskip It is well known that the tube domain $\6H$ has a canonical
realization as bounded symmetric domain $$\6D:=\{z\in
E:\One_{n}-z\overline z\hbox{ positive definite}\}$$ given by the
Cayley transformation $\gamma:z\mapsto
(z-\One_{n})(z+\One_{n})^{-1}$. The boundary $\partial\6D$ is the
union of the CR-submanifolds of $E$ $$\6R_{p,0}:=\{z\in
E:(\One_{n}-z\overline z)\in\6C_{p,0}\},\qquad0\le p<n\,.$$ The group
$\Aut(\6D)$ acts transitively on each $\6R_{p,0}$ and can be
identified with $\Aut(\6R_{p,0})$ this way for every $p>0$, see
\Lit{KAZT}. The Cayley transformation $\gamma$ gives a CR-isomorphism
from $\6T_{p,0}$ to the dense domain
$\{z\in\6R_{p,0}:\det(\One_{n}-z\overline z)\ne0\}$ in $\6R_{p,0}$.
Furthermore, the manifold $\6R_{p,0}$ is an $\Aut(\6H)$-equivariant
fiber bundle over $\U(n)/\big(\O(n{-}p)\times\U(p)\big)$ with fibers
the holomorphic arc components. In particular,
$\pi_{1}(\6R_{p,0})=\ZZ_{2}$ for $p<n$. As in the proof of Proposition
\ruf{BA} it is shown for every $p>0$ that every continuous CR-function
on the universal covering $\mu:\tilde\6R_{p,0}\to\6R_{p,0}$ is
constant on $\mu$-fibers.

Every $\6T_{p,0}$ has the following local description in terms of
normalized coordinates generalizing \Ruf{AO}: For $k:=n-p$ put
$$W:=\{w\in\CC^{k\Times k}:w'=w\}\,,\qquad E_{1}:=\CC^{p\times
k}\,,\qquad E_{2}:=\{z_{2}\in\CC^{p\Times p}:z_{2}'=z_{2}\}$$ and
denote by $\6M$ in $U:=\{(w,z_{1},z_{2})\in W\times E_{1}\times
E_{2}:\One_{p}-z_{2}\overline z_{2}>0\}$ the CR-submanifold given by
the matrix equation $$w+\overline w\;=\;x+\overline
x\Steil{with}x:=\overline z_{1}'(\One_{p}-z_{2}\overline
z_{2})^{-1}(z_{1}+z_{2}\overline z_{1})\,.\Leqno{OA}$$ Then
$$(w,z_{1},z_{2})\;\longmapsto\;\pmatrix{
w+z_{1}'yz_{1}&\sqrt2z_{1}'y\cr \noalign{\vskip4pt}
\sqrt2yz_{1}&(\One_{p}-z_{2})y\cr
}\Steil{with}y:=(\One_{p}+z_{2})^{-1}\Leqno{OB}$$ defines a
biholomorphic mapping $\phi$ from $U$ to a domain in $E$ such that
$\phi(\6M)$ is open in $\6T_{p,0}$. This follows from the realization of
$\6H$ as Siegel domain of the third kind in terms of the partial
Cayley transformation \Ruf{OB}, compare \Lit{PIAT} and \Lit{LOSO}
\p10.7.

\bigskip Due to the coincidence $\7{sp}(2,\RR)\cong\7{so}(2,3)$
between symplectic and orthogonal groups in low dimensions there are
two ways to generalize the tube $\6T$ over the light cone
$\6C\subset\RR^{3}$ to higher dimensions: The first possibility in
terms of symmetric $n\Times n$-matrices and the symplectic group has
been explicitly described above. The corresponding tube manifolds
$\6T_{p,q}$, $0<p+q<n$, are all not simply connected and hence provide
nontrivial examples for Proposition \ruf{OD}. Let us close with a
description of the second way associated with orthogonal groups: For
fixed $n\ge2$ let $$\langle z|w\rangle:=z_{1}w_{1}+\dots+z_{n}w_{n}$$
be the standard symmetric bilinear form and $z\mapsto\overline
z=(\overline z_{1},\dots,\overline z_{n})$ the standard conjugation on
$\CC^{n}$. Denote by
$$\6C:=\big\{(t,x)\in\RR^{2}\Times\RR^{n}:t_{1}t_{2}=\langle
x|x\rangle\,,\;t_{1}\!+t_{2}>0\big\}$$ the future light cone in
$(n+2)$-dimensional space time and by $\6T:=\6C\oplus i\RR^{n+2}$ the
corresponding tube manifold over $\6C$. Then $\6T$ is a simply
connected, $2$-nondegenerate homogeneous CR-manifold with
$\hol(\6T,a)\cong\7{so}(2,n{+}2)$ for all $a\in\6T$, see \Lit{KAZT}
for the last statement. Then the hypersurface $\6M$, given in
$U:=\{(w,z_{1},z_{2})\in\CC\Times\CC^{n}\Times\CC:|z_{2}|<1\}$ by the
equation $$w+\overline w=\big(2\langle z_{1}|\overline
z_{1}\rangle+\langle z_{1}|z_{1}\rangle\overline
z_{2}+\langle\overline z_{1}|\overline z_{1}\rangle
z_{2}\big)(1-z_{2}\overline z_{2})^{-1}\,,$$ is CR-equivalent to a
domain in $\6T$. As in \Ruf{AO} an explicit equivalence is obtained by
$$(w,z_{1},z_{2})\mapsto(1+z_{2})^{-1}\big(w+wz_{2}+ \langle
z_{1}|z_{1}\rangle,\,1-z_{2},\,\sqrt2z_{1}\big)\,.$$

\medskip

\noindent{\bf Acknowledgment:} The authors wish to thank D. Zaitsev for
helpful discussions.

{\parindent 15pt\bigskip\bigskip\Klein {\noindent\gross References}
\bigskip

\font\klCC=txmia scaled 900

\Ref{AZAD}Azad, H., Huckleberry, A., Richthofer, W.: Homogeneous CR-manifolds. J. Reine Angew. Math. {\bf 358} (1985), 125--154.
\Ref{BERO}Baouendi, M.S., Ebenfelt, P., Rothschild, L.P.: {\sl Real Submanifolds in Complex Spaces and Their Mappings}. Princeton Math. Series {\bf 47}, Princeton Univ. Press, 1998.
\Ref{BGSS}Boggess,~A.: {\sl CR Manifolds and the Tangential Cauchy-Riemann Complex}. Studies in Advanced Mathematics. CRC Press. Boca Raton~ Ann Arbor~ Boston~ London 1991.
\Ref{BUSH}Burns, D., Shnider, S.: Spherical hypersurfaces in complex manifolds. Invent. Math. {\bf 33} (1976), 223--246. 
\Ref{CHMO}Chern, S.S., Moser, J.K.: Real hypersurfaces in complex manifolds. Acta. Math. {\bf 133} (1974), 219-271.
\Ref{EBFT}Ebenfelt, P.: Normal Forms and Biholomorphic Equivalence of Real Hypersurfaces in \hbox{\klCC\char131}$^3$. Indiana J. Math. {\bf 47} (1998), 311-366.
\Ref{EBEN}Ebenfelt, P.: Uniformly Levi degenerate CR manifolds: the 5-dimensional case. Duke Math. J. {\bf 110} (2001), 37--80.
\Ref{ERSH}Ershova, A.E.: Automorphisms of $2$-Nondegenerate Hypersurfaces in \hbox{\klCC\char131}$^3$. Math. notes {\bf 69} (2001), 188-195.
\Ref{GAME}Gaussier, H., Merker, J.: A new example of a uniformly Levi degenerate hypersurface in \hbox{\klCC\char131}$^3$. Ark. Mat. {\bf 41} (2003), 85--94.
\Ref{KAPP}Kaup, W.: On the holomorphic structure of $G$-orbits in compact hermitian symmetric spaces. Math. Z. {\bf 249} (2005), 797-816.
\Ref{KAZT}Kaup, W., Zaitsev, D.: On the local CR-structure of Levi-degenerate group orbits in compact Hermitian symmetric spaces. J. Eur. Math. Soc., to appear 
\Ref{KLIN}Klingen, H.: Diskontinuierliche Gruppen in symmetrischen R\"aumen. Math. Ann. {\bf 129} (1955), 345-369.
\Ref{LOSO}Loos, O.: {\sl Bounded symmetric domains and Jordan pairs.} Mathematical Lectures. Irvine: University of California at Irvine 1977. 
\Ref{PALA}Palais, R.S.: {\sl A global formulation of the Lie theory of transformation groups.} Mem. AMS 1957.
\Ref{PIAT}Piatetsky-Chapiro, I.I.: {\sl G\'eom\'etrie des domaines classiques et th\'eorie des fonctions automorphes.} Dunod, Paris 1966.
\Ref{WOLF}Wolf, A.J.: {\sl Fine Structure of Hermitian Symmetric Spaces.} Symmetric spaces (Short Courses, Washington Univ., St. Louis, Mo., 1969--1970), pp. 271--357. Pure and App. Math., Vol. 8, Dekker, New York, 1972.
\par}

\bigskip\bigskip\noindent Mathematisches Institut, Universit\"at
T\"ubingen, Auf der Morgenstelle 10, 72076 T\"ubingen, Germany
(e-mails: {\tt kaup@uni-tuebingen.de}, {\tt fels@uni-tuebingen.de})

\closeout\aux\bye